\documentclass[twoside,11pt]{article}
\usepackage{jmlr2e}
\usepackage{xcolor}
\usepackage{microtype}
\usepackage{graphicx}
\usepackage{subfigure}
\usepackage{booktabs} 

\usepackage{mathrsfs,amsmath, amscd,thmtools}
\usepackage{stmaryrd}
\usepackage{arydshln}
\DeclareMathOperator*{\argmax}{arg\,max}
\DeclareMathOperator*{\argmin}{arg\,min}
\usepackage{blindtext}
\usepackage{hyperref}
\usepackage{algorithm}
\usepackage[noend]{algpseudocode}
\usepackage{pifont}
\usepackage{sidecap}
\usepackage{wrapfig}
\usepackage{enumerate}
\usepackage{enumitem}
\usepackage{comment}
\usepackage{blindtext}
\usepackage{hyperref}
\usepackage[mathscr]{euscript}

\newcommand{\Trace}{\operatorname{Tr}}

\newcommand\SL{\text{SL}}
\newcommand\GL{\text{GL}}
\newcommand\R{\mathbb{R}}

\def\be{\begin{equation}}
\def\ee{\end{equation}}
\def\bes{\begin{eqnarray}}
\def\ees{\end{eqnarray}}

\def\mc{\tau}

\def\2{\frac{1}{2}}
\def\4{\frac{1}{4}}

\newcommand{\T}{\mathcal}
\newcommand{\bU}{\mathbf{U}}
\newcommand{\Lc}{\mathbf{\bar{\mathcal{L}}}}
\newcommand{\m}[1]{{\bf{#1}}}
\newcommand{\wh}[1]{{\widehat{#1}}}

\renewcommand{\mc}[1]{\ensuremath{\mathcal{#1}}} 

\newcommand{\mb}[1]{{\mathbb{#1}}}

\firstpageno{1}

\begin{document}

\title{A Bilevel Optimization Method for Tensor Recovery Under Metric Learning Constraints}

\author{\name Maryam Bagherian \email bmaryam@umich.edu \\
       \addr Department of Computational Medicine and Bioinformatics\\
       \addr Michigan Institute for Data Science (MIDAS)\\
       University of Michigan, Ann Arbor, MI 48109, USA
       \AND
       \name Davoud A. Tarzanagh \email tarzanaq@umich.edu \\
       \addr Department of Electrical Engineering and Computer Science\\
       University of Michigan, Ann Arbor, MI 48109, USA
   \AND
   \name Ivo Dinov \email statistics@umich.edu \\
   \addr Department of Health Behavior and Biological Sciences \\
   Michigan Institute for Data Science (MIDAS)\\
   University of Michigan, Ann Arbor, MI 48109, USA
\AND
\name Joshua D. Welch \email welchjd@umich.edu \\
\addr Department of Computational Medicine and Bioinformatics\\
Department of Electrical Engineering and Computer Science\\
University of Michigan, Ann Arbor, MI 48109, USA}

\editor{}

\maketitle

\begin{abstract}
Tensor completion and tensor decomposition are important problems in many domains. In this work, we leverage the connection between these problems to learn a distance metric that improves both decomposition and completion. We show that the optimal Mahalanobis distance metric for the completion task is closely related to the Tucker decomposition of the completed tensor. Then, we formulate a bilevel optimization problem to perform joint tensor completion and decomposition, subject to metric learning constraints. The metric learning constraints also allow us to flexibly incorporate similarity side information and coupled matrices, when available, into the tensor recovery process. We derive an algorithm to solve the bilevel optimization problem and prove its global convergence. When evaluated on real data, our approach performs significantly better compared to previous methods.
\end{abstract}


\section{Introduction}\label{sec:intro}

Tensor factorization and tensor completion have gained increased attention and popularity over the last decade. These methods have found broad applications in many fields  \citep{kolda2009tensor,cichocki2015tensor,bi2020tensors}. The tensor completion problem seeks to find a tensor that is simultaneously close to an incompletely observed tensor and parsimonious, in the sense that the completed tensor has low rank~\citep{kolda2009tensor,cichocki2015tensor,bi2020tensors}. This is directly analogous to the matrix completion problem, for which there are many effective and robust methods \citep{davenport2016overview}. There is a close connection between tensor completion and tensor decomposition, which has emerged as a key tool for exploratory analysis of high-dimensional datasets. Popular types of tensor decomposition include Tucker \citep{tucker64extension}, Canonical Polyadic (CP) \citep{carroll1970analysis}, higher-order SVD (HOSVD) \citep{de2000multilinear}, tensor train (TT) \citep{oseledets2011tensor}, and tensor SVD (t-SVD) \citep{kilmer2013third}. 

Generalizing matrix factorization techniques, tensor decomposition represents a tensor as a product of more basic and interpretable components. In many scientific applications, a key objective is to find components with physical interpretations that represent latent data-generating variables. However, there are several key challenges related to decomposing multi-dimensional multi-modal datasets: (i)  missing values, (ii) the need to find sparse representations, and (iii) heterogeneity.\\ 

\noindent (i). \textbf{Tensor completion} methods were developed to address the issue of missing values \citep{liu2012tensor,kressner2014low,zhao2015bayesian,zhang2014novel,song2017tensor, bagherian2021coupled}. Most existing tensor completion
work deals with sampling procedures that are independent
of the underlying data values. While this assumption allows
the derivation of nice theoretical guarantees, it seldom holds
in real-world applications. When some \emph{regularity} information is present, one may utilize it in order to better predict the missing entries. The recent method proposed in \citep{bagherian2021coupled} is a pertinent example where additional information are available and used for a better prediction. As another example, one may refer to \citep{narita2012tensor} where authors proposed two regularization methods called ``within-mode regularization" and ``cross-mode regularization" to incorporate auxiliary regularity information in the tensor completion problems. The key idea is to construct within-mode or cross-mode regularity matrices, incorporate them as smooth regularizers, and then combine them with a Tucker decomposition to solve the tensor completion problem.\\

\noindent (ii). \textbf{Sparse representation} aims to solve high-dimensional tensor problems, and shrinkage the size of tensor models. Tucker and Kronecker models of tensors allow sparse representation of a tensor and can be used in order to better model and interpret such data. A number of Kronecker-based dictionary learning methods have been proposed in literature, including   \citep{hawe2013separable,qi2018multi,shakeri2018minimax,bahri2018robust} and associated algorithms. \\

\noindent (iii). \textbf{Heterogeneity of tensor models} is perhaps the challenge least discussed in the literature. However, in a number of applications, data heterogeneity is a key property that must be modeled. For example, image data exhibit heterogeneity due to differences in lighting and posing,
which can be extracted a priori from available metadata and utilized at analysis time  \citep{tarzanagh2019regularized}.

This paper addresses these challenges and
develops a bilevel optimization method for tensor recovery under metric learning constraints. Our proposed approach combines regularized Tucker factorization with a nonconvex tensor model that allows for high-dimensional tensor completion, and metric learning constraints to improve the quality of tensor recovery in the heterogeneous setting. Specifically, we make multiple algorithmic and theoretical contributions summarized below.

\begin{itemize}
	\item
	[I.] We draw new connections between tensor completion and tensor decomposition, showing that the two problems can be related by a common Mahalanobis distance metric. 

    \item [II.] We formulate a bilevel optimization for tensor recovery problems, called Tensor Recovery Under Metric Learning Constraints (TR-MLC), which uses a distance metric learned jointly from tensor completion and decomposition to improve both. In the proposed bilevel problem, the upper level problem uses regularized Tucker model to allow sparse representation of the underlying high-dimensional data, and the lower level problem uses a set of coupled matrices to inform the tensor completion task. This allows simultaneous completion and decomposition of a tensor while predicting missing entries.
	\item [III.] We develop an alternating direction method of multipliers (ADMM) algorithm for solving the bilevel optimization problem and prove its global convergence.
    \item [IV.] We demonstrate the advantages of TR-MLC compared to existing approaches on real datasets.
\end{itemize}

A key motivation for the proposed method is that both tensor completion and tensor decomposition are important in many domains, including chemometrics \citep{smilde2005multi}, audio and speech processing \citep{wu2010robust}, neuroscience \citep{beckmann2005tensorial}, and computational biology \citep{bagherian2021coupled}. In such applications, both completion and decomposition of the same tensor are often of interest. 

The rest of the manuscript is organized as follows: Related literature is discussed in Section~\ref{sec:RL} and a brief background on the methodology is provided in Section~\ref{sec:BG}. We provide the model framework in Section~\ref{sec:MF} followed by an ADMM algorithm in Section~\ref{sec:ADMM}, and we discuss the convergence of the proposed model in Section~\ref{sec:globconv}. Finally, we test the performance of TR-MLC on a number of real datasets in Section~\ref{sec:ER} and 
we conclude in Section~\ref{Conclusion}. Additional details, when necessary, are provided in the Appendix. 
\section{Related Literature}\label{sec:RL}

\textbf{Similarity-based learning algorithms} are among the first methods to be generalized for tensors in the field of machine learning \citep{yan2005discriminant,fu2008image}. Analogous to those for matrices, they are designed based on the algorithms that can learn from a dataset according to the similarities present between the data. Choosing an appropriate distance learned from the data has proven to be able to greatly improve the results of
distance-based algorithms \citep{yang2006distance}. A good distance allows data to be transformed to facilitate their analysis, with mechanisms such as dimensionality reduction and/or completion \citep{suarez2018tutorial}.

\noindent \textbf{Graph-Regularized Tensor Methods.}
Numerous previous papers have incorporated graph regularization terms into either the tensor completion or tensor decomposition objective. For instance, authors in \citep{takeuchi2016graph} introduced a graph Laplacian based regularizer and used it for inducing latent factors to represent auxiliary structures. Li et al. used a graph-regularized tensor completion model to impute spatial transcriptomic data \citep{li2021imputation}. The model is regularized by a product of two chain graphs which intertwines the two modes of the tensor. Many additional related papers have used graph-based side information to regularize tensor completion or decomposition  \citep{guan2020alternating,song2021dstg,sofuoglu2020graph}.

\noindent\textbf{Distance Metric Learning and Tensors.}
Learning a distance metric in a feature space dates back to one the simplest methods for classification in 1951 and emerged in other approaches such as support vector machines and K-means clustering. However, these types of algorithms have widely used standard metrics, such as the Euclidean metric, which may not accurately capture the important relationships within the data \citep{weinberger2009distance,yang2006distance}. Distance metric learning was originally developed in the context of supervised problems such as classification. Metrics are often learned from pairs of points explicitly labeled as similar or dissimilar \citep{xing2002distance}. The technique has been applied to tensors in several studies, such as \citep{hao2019collect,hou2007saliency,cai2005subspace}. However, these papers were also focused primarily on supervised learning, rather than tensor decomposition or completion. 

\noindent \textbf{Simultaneous Tensor Decomposition and Completion.}
The approach most closely related to ours is Simultaneous Tensor Decomposition and Completion (STDC) \citep{chen2013simultaneous}. STDC completes a partially observed tensor by finding a Tucker decomposition that exactly recovers the known entries. STDC further incorporates graph-structured side information using a graph regularization penalty to constrain the factor matrices of the Tucker decomposition. Our approach is related in the sense that we also perform joint decomposition and completion, and our lower-level optimization problem contains a trace term similar to the graph regularization penalty. However, a key difference is that our lower-level objective function includes a distance metric learning term, which significantly changes the completion objective. As described in the results section, when we augment the STDC objective with a similar term, it significantly improves the performance of STDC. Furthermore, we formulate our approach as a bilevel optimization problem, which lays a flexible foundation for future extensions. This allows us to, for example, incorporate coupled matrices that are not either graphs or similarity matrices.

\section{Background}\label{sec:BG}
\subsection{Tensor Completion}

The goal of tensor completion is to estimate an unknown tensor $\wh{\T{X}} \in \mathbb R^{N_1\times N_2\times \cdots \times N_K}$ from an incomplete tensor $\T X$ where the set of indices for the observed entries is $\m \Omega_{\T X}$. Throughout, we aim to find a \textit{low-rank} approximation $\T{\wh X}$ that perfectly recovers the observed entries of $\T X$. 
A  natural approach to completing tensors is to solve the following convex optimization problem: 
\begin{equation}
\label{eq:1}
\begin{aligned}
\min\limits_{ \T{\wh X} }\quad \|\T{\wh X}\|_{*} \quad 
\textrm{s.t.} \quad \m \Omega_{\T X} = \m \Omega_{\T{\wh X}}.
\end{aligned}
\end{equation}
where $ \|\T{\wh X}\|_{*}$ is the tensor nuclear norm. Unfortunately, the tensor nuclear norm is hard to compute so this approach does not lead to any algorithmic guarantees \citep{hillar2013most}.
\subsection{Tensor Decomposition}
Tucker decomposition \citep{tucker1964extension}, decomposes a tensor into a set of matrices and a smaller core tensor. In the following, we provide a formal definition of this problem.
\begin{definition}[Tucker Decomposition]\label{def:dcot}
Given an $K$-way tensor $\T Z\in \mathbb R^{N_1\times N_2\times \cdots \times N_K}$, its \textnormal{Tucker} decomposition is given by
\begin{align}\label{eq:dcot}
\T{Z} &= \T{G}\times_1 \m {V}^{(1)}\times_2 \cdots \times_K \m {V}^{(K)},
\end{align}
where $\m{V}^{(\ell)} =[\m{v}_1^{(\ell)},\m{v}_2^{(\ell)},\cdots,\m{v}_{n_\ell}^{(\ell)}] \in \R^{N_n\times n_\ell}$, is the $\ell$-th factor matrix consisting of latent components $\m{v}_k^{(\ell)}$ and $\T{G} \in \R^{n_1 \times n_2 \times \cdots \times n_K}$ is a  core tensor reflecting the connections (or links) between the latent components and factor matrices.
\end{definition}

In Definition~\ref{def:dcot}, the $N$-tuple $(n_1, n_2, \dots , n_N)$ with $n_n = \text{rank}(\m{Z}_{(n)})$ is called the multi-linear rank of $\T{Z}$. For a core tensor of minimal size, $n_1$ is the column rank (the dimension of the subspace spanned by mode-1 fibers), $n_2$ is the row rank (the dimension of the subspace spanned by mode-2 fibers), and so on. An important difference from the matrix case is that the values of $n_1,n_2,\cdots,n_N$ can be different for $N\geq3$. Note that decomposition \eqref{eq:dcot} can be expressed in a matrix form as:
\begin{eqnarray}\label{eq:dcotmode}
 \m{Z}_{(\ell)} & =&  \bU^{(\ell)} (\m{G}_{(\ell)}) (\bigotimes_{k \neq \ell} \m{U}^{(k)})^\top.
\end{eqnarray}
The Tucker model formulation provides a generic tensor decomposition that encompasses many other popular tensor decomposition models. Indeed, when $\T{H}=0$ and $\m{V}^{(\ell)}$ for $ \ell \in \{ k=1, 2, \cdots, K\}$ are orthogonal, \eqref{eq:dcot} corresponds to HOSVD \citep{de2000multilinear}. The CP decomposition \citep{carroll1970analysis} can also be considered as a special case of the Tucker model with super-diagonal core tensors.

\subsection{Distance Metric Learning}
Distance metric learning \citep{suarez2018tutorial} aims to learn distances from the data, where distance refers to a map $d:A\times A\to \mb R_+$, {$A$ is a non-empty set and $d$ satisfies the following conditions:}

\begin{enumerate}
	[label=(\roman*)]
	\item \label{dis:co}   Coincidence: $d(\m a,\m b)=0~\iff~\m a=\m b, \forall \m{a}, \m{b} \in A.$ 
	\item \label{dis:sym}  Symmetry: $d(\m a,\m b)=d(\m b,\m a), ~  \forall \m{a}, \m{b} \in A.$
	\item \label{dis:tri}  Triangle inequality: $d(\m a,\m c)\le d(\m a,\m b)+d(\m b,\m c),  \forall \m{a}, \m{b}, \m{c} \in A.$
\end{enumerate}
Other properties such as \emph{non-negativity}, \emph{reverse triangle inequality}, and \emph{generalized triangle inequality} follow immediately from the definition above. 
Distance metric learning frequently focuses on learning Mahalanobis distances, since they
are parametrized by matrices, and therefore are computationally tractable. Mahalanobis distances satisfy additional properties, including \emph{translation invariance} and \emph{homogeneousness} \citep{suarez2018tutorial}. In a $K$-dimensional Euclidean space, we may form a family of metrics over the set $A$ by computing Euclidean distances after a linear transformation $\m a\to\m L(\m a)$, where $\m L$ is injective. Therefore, the squared distances can then be computed as $d(\m a,\m b)=\|\m L(\m a-\m b)\|^2_F$ which may also be expressed in terms of a positive semidefinite square
matrix $\m M = \mathbf L\mathbf L^\top$. If $\m L$ is also surjective, which results in $\m M$ being full rank, the matrix $\m M$ parametrizes the distance $d$. Then the matrix $\m M$ is referred to as a \emph{Mahalanobis} metric. In Gaussian distributions, the matrix $\m M$ plays the role of the inverse covariance matrix.

\begin{definition}[Tensor Mahalanobis Distance]\label{def:maha_t}
	Let $\T X\in \mathbb R^{N_1\times N_2\times \cdots\times N_K}$. Consider the multilinear transformation\\ $ \varphi:\mathbb R^{N_1\times N_2\times \cdots\times N_K}\to \mathbb R^{N_1\times N_2\times \cdots\times N_K}$, with
	\begin{align}
	\varphi(\T X)= \T X\times_1 \m L^{(1)}\times_2 \m L^{(2)}\cdots \times_K \m L^{(K)}, 
	\end{align}
	where the square matrices $\m L^{(\ell)}\in \mathbb R^{N_\ell\times N_\ell}$, for $\ell=1,\cdots, K$, are called the $\ell$-mode matrices. The squared Mahalanobis distance can be computed as:
	\begin{align}
	\label{eq:norm:def}
	\nonumber 
	d_M(\T X_i, \T X_j) &:= \big{\|} \varphi(\T X_i)- \varphi(\T X_j)\big{\|}_F^2\\
	&= \big{\|}(\T X_i -\T X_j)\times_1 \m L^{(1)} \cdots \times_K \m L^{(K)}\big{\|}_F^2.
	\end{align}
\end{definition}

\begin{lemma}[\cite{suarez2018tutorial}]
	{Let  $\T X_i, \T X_j \in \mathbb R^{N_1\times N_2\times \cdots\times N_K}$. Then} 
	\begin{equation}\label{eq:maha:def}
	d_M(\T X_i, \T X_j)=
	\Trace\Big(\wh{\m L}^{(\ell)}(\m X_i-\m X_j)_{(\ell)} \, 
	\wh{\T L}^{\neq \ell}_{\otimes}(\m X_i-\m X_j)_{(\ell)}^\top\Big ),
	\end{equation}
	where $\m X_{(\ell)}$ denotes the $\ell$-th matricization of tensor $\T X$, for each $\ell$, $\wh{\m L}^{(\ell)}= \m L^{(\ell)^\top}\m L^{(\ell)},$ and $\wh{\T L}_\otimes^{\neq \ell}=\wh{\m L}^{(K)}\otimes \cdots \times \wh{\m L}^{(\ell-1)}\otimes \wh{\m L}^{(\ell+1)}\otimes \cdots\otimes \wh{\m L}^{(1)},$ where $\otimes$ denotes Kronecker product. 	Here, if $\m L^{(\ell)}$, for $\ell=1, \cdots, K$ are orthogonal matrices, then $ d_M$ recovers Euclidean distance. 
\end{lemma}


\section{Tensor Recovery Under Metric 
	Learning Constraints}\label{sec:MF}

In this section, we introduce TR-MLC which uses a distance metric learned jointly from tensor completion and decomposition.

\subsection{Tensor Completion with Mahalanobis Metric}
We formulate the problem (introduced in Eq.~\eqref{eq:1}) of completing a low rank tensor in the sense of a Mahalanobis metric defined in Eq.~\eqref{eq:norm:def}. That is, we want to find a completed tensor $\hat{\T X}$ that exactly matches the original tensor in the observed entries and minimizes $d_M(\hat{\T X}, \mathbf{0})$, {where $\mathbf{0}$ denotes the tensor of all zeros}. This gives the following optimization problem:   
\begin{equation}\label{eq:main:1}
\begin{aligned}
\min\limits_{ \Gamma} \quad &  \frac12\, \big{\|}\T {\wh X}\times_1 \m { L}^{(1)}\times_2 \m { L}^{(2)}\cdots \times_K \m {L}^{(K)}\big{\|}_F^2 \\
\textrm{s.t.} \quad&   \m \Omega_{\T X}= \m \Omega_{\T { \wh X}}, 
\end{aligned}
\end{equation}
where the collection of square matrices $\{\m { L}^{(\ell)}\}_{\ell=1}^K$, $\m { L}\in \mathbb R^{N_\ell \times N_\ell}$, form a family of linear transformations as in Definition~\ref{def:maha_t} and $\Gamma=\big\{\T { \wh X},\m L^{(1)}, \m L^{(2)}, \cdots, \m L^{(K)}\big\}$, for all $\ell \in \{1, \cdots, K\}$. 

If side information about the similarity of each tensor mode is available, we can also incorporate it to further refine the metric by adding a regularization term to the objective. Suppose that for any fixed mode $\ell \in \{1,\cdots, K\}$, there exists a square symmetric matrix, $\m S_{(\ell)}\in \mathbb R^{N_\ell \times N_\ell}$, containing the similarity information regarding the specific mode $N_\ell$.  {Without loss of generality}, we may assume that $\m S_{(\ell)}$'s
are positive semidefinite. The completion Problem \eqref{eq:main:1} can then be rewritten to incorporate the similarity information as follows:
\begin{equation}\label{eq:main:2}
\begin{aligned}
\min\limits_{\Gamma  } \quad  & \frac12 \, \big{\|}\T {\wh X}\times_1 \m L^{(1)}\times_2 \m L^{(2)}\cdots \times_K \m L^{(K)}\big{\|}_F^2\\ 
\quad  &+\frac12\,  \sum\limits_{\ell=1}^K \lambda_\ell\, \Trace(\m L^{(\ell)}\, \m S_{(\ell)}\, \m L^{{(\ell)^\top}}), \\
\textrm{s.t.} \quad  & \m \Omega_{\T X}= \m \Omega_{\T {  \wh X}}. 
\end{aligned}
\end{equation}
In many data analysis applications where tensor decomposition is extensively used \citep{liu2012tensor,chen2013simultaneous,kolda2009tensor,tarzanagh2019regularized}, the multi-dimensional data exhibit sparse, group- sparse or other type of specialized structures. Such properties are related to the sparse representation and heterogeneity challenges that we noted in the introduction. This motivates the use of \textit{penalized} tensor decomposition, where some type of nonsmooth regularization such as $\ell_1$-norm, $\ell_2$-norm, and nuclear norm is imposed on either core or factor matrices.

In the light of this observation, we next develop a bilevel optimization problem that links penalized Tucker factorization with tensor completion.

\subsection{Mahalanobis Metric Links Tensor Completion and Decomposition}\label{subsec:MD-TC}
A key observation about the optimization problem \eqref{eq:main:2} led us to initially pursue a joint completion and decomposition approach. That is, the linear transformation defining a Mahalanobis metric is itself a Tucker decomposition. To see this, note the similarity between Eqs.~\eqref{eq:dcot} and  ~\eqref{eq:norm:def}. The matrices $\m L^{(\ell)}$ defining a Mahalanobis metric can thus be interpreted as the matrices $\m V^{(\ell)}$ of a Tucker decomposition. We further had the intuition that completing a tensor should lead to more accurate decomposition of the tensor into underlying factors. These insights led us to formulate a bilevel optimization problem in which decomposition is performed jointly with completion, and both problems are linked by a learned distance metric.

Another motivation for our formulation is a connection between Mahalanobis distance and group theory. Recall that Mahalanobis distance is equivalent to Euclidean distance after a linear change of coordinates that normalizes the covariance matrix/tensor of the data to the identity. The Kempf-Ness theorem \citep{olver1999classical}, a classic result from group theory, guarantees that the problem of finding such a change of coordinates has a unique solution that is optimal in a certain sense. Briefly, a proper action on an arbitrary group $H$ can perform the change of coordinates in the vector space $V$ such that it preserves the mathematical structure of the elements. It is known that the group $H$ has a unique maximal compact subgroup $ M$. 
\begin{theorem}[Kempf-Ness Theorem]\label{thm:KNT}
	Consider the real-valued map ~$\varphi: H\to \R$ given by ~$\varphi(h)=\|h\cdot \mathbf v\|^2$ for some nonzero  $ \mathbf v\in \R^n$. If a critical point exists for $\varphi$, then it is an absolute minimum and $ Mh$ is a coset of group $ H$~ for $h\in H$. Moreover, the set of critical points is unique up to the action of group $M$.
\end{theorem}
Further explanation of this result and connections to the representation theory of reductive groups can be found in the Appendix, Section~\ref{app:kn}. 

These insights led us to formulate a bilevel optimization problem in which decomposition is performed jointly with completion, and both problems are linked by a learned Mahalanobis distance metric.

\subsection{Bilevel Optimization}\label{subsec:bo}
 The technique of bilevel optimization has gained popularity in recent years. Bilevel optimization problem are often used to model hierarchical processes in which an upper level problem makes a decision anticipating the rational inputs of the lower level problem. In general, any bilevel optimization is locally equivalent to a single-level optimization problem \citep{dempe2014new}. Bard \citep{bard2013practical} studied the linear bilevel problem
and developed first-order necessary optimality conditions under which any linear bilevel optimization problem can be formulated as an equivalent single-level linear optimization problem.

The most direct approach to solving a bilevel optimization problem is to solve the equivalent mathematical problem given in Eq~\eqref{prob:1} (see Appendix~\ref{app:bl}) .  For any linear bilevel problem, a single level problem exists where the solutions of both problems are equivalent \citep[ch.5]{bard2013practical}. Further details for a general linear problem is provided in Appendix~\ref{app:bl}. In what follows we propose a bilevel optimization problem followed by a linearized ADMM approach to solve the proposed problem.  

\subsection{A Bilevel Optimization Framework for Tensor Recovery Under
Metric Learning Constraints}\label{subsec:OP}
Let $\mathbf{\mc{X}}\in \mathbb R^{N_1\times \cdots \times N_K}$ be a tensor that admits a Tucker factorization. Our goal is to perform tensor recovery by jointly finding a completed tensor $\wh{\T X}$, a penalized Tucker factorization {$\llbracket \T G; \{\m{V}^{(\ell)}\}_{\ell=1}^K \rrbracket$} of $\mathbf{\mc{X}}$, and Mahalanobis metric defined by $\{\m{L^{(\ell)}}\}_{\ell=1}^K$. We formulate a penalized tensor factorization under metric learning constraints as the following bilevel optimization problem:

\begin{subequations}\label{eq:bilTen}
	\begin{align}
	\label{eq:bil:up:Ten}
	\min\limits_{ \T{Z} } &~~ 
	\frac{1}{2} \| \T X- \T{Z}\|^2_\textnormal{F}+ J_1(\T{G}) + \sum\limits_{\ell=1}^K J_{(2,\ell)} ({\m{V}}^{(\ell)}) \\
	\label{eq:bil:low:Ten}
	\textrm{s.t.}~~&~ \T{Z}  :=~\argmin\limits_{\wh{\T{Z}}} ~  \frac{1}{2}\| \wh{\T{ X}}\times_1 \hat{\m L}^{(1)}\times_2\cdots\times_K \hat{\m L}^{(K)}\, \|^2_\textnormal{F}\\
	\nonumber
	&\qquad  \qquad \qquad+ \sum\limits_{\ell=1}^K\lambda_{\ell}
	\Trace (\hat{\m L}^{(\ell)}\,  \m{S}_{(\ell)} \, \m{\hat L}^{(\ell)^\top}),
	\end{align}
	and 
	\begin{equation}\label{eq:prob:form}
	\begin{cases}
	\m\Omega_{\T X}= \m\Omega_{\wh{\T X}},\\ 
	\m{\hat{L}^{(\ell)}}=
	\m{\hat V}^{(\ell)^\top} \, \m{\hat V}^{(\ell)},\\
	\wh{\T G} = \wh{\T X}\times_1
	\m {\hat V}^{(1)^\top} \times_2 \cdots \times_K 
	\m {\hat V}^{(K)^\top}, \\
	\wh{\T Z} = \wh{\T G}\times_1 \hat{\m V}^{(1)}\times_2\cdots \times_K \hat{\m V}^{(K)}, \\
	{\T Z} = {\T G}\times_1 {\m V}^{(1)}\times_2\cdots \times_K {\m V}^{(K)}.
	\end{cases}
	\end{equation}
\end{subequations}
We detail the purpose of each component of \eqref{eq:bilTen}:

$\bullet$ The lower-level problem \eqref{eq:bil:low:Ten} corresponds to low-rank tensor completion introduced in \eqref{eq:main:2} which allows incorporating side information about the similarity of each tensor model. In \eqref{eq:bil:low:Ten}, the set $\{\lambda_\ell\}_{\ell=1}^K$ denotes the lower-level regularization parameters. Note that if $\{\m S_{(\ell)}\}_{\ell=1}^K$ are available, the second term in \eqref{eq:bil:low:Ten} uses them to refine the distance metric. Even if the $\{\m S_{(\ell)}\}_{\ell=1}^K$ matrices are not given, we can initialize them by, for instance, taking pairwise distances along each mode of the tensor \citep{narita2012tensor}. We can then iteratively update them using pairwise distances from the completed tensor at each stage of optimization process.

$\bullet$ The upper-level problem \eqref{eq:bil:up:Ten} is penalized Tucker decomposition to solve high-dimensional tensor problems. We use $J_1(\T{G})$ and {$J_{(2,\ell)} ({\m{V}}^{(\ell)})$} to denote penalty functions on $\T{G}$ and {$\m{V}^{(\ell)}$}, respectively.  Note that the core tensor is the projection of  $\T X$ onto the tensor basis formed by the factor matrices ${\m V}^{(\ell)}$; that is, ${\T G}=\T X \bigotimes_{\ell=1}^K \m V^{(\ell)^\top}$. If the core tensor $\T G$ has the same dimensions as the data tensor $\T X$, we can simply set $\m{\hat L}^{(\ell)}=\m{\hat V}^{(\ell)}$. In this case, $\{\m{L^{(\ell)}}\}_{\ell=1}^K$ define a true Mahalanobis metric. However, the above formulation also allows us to learn a core $\T G$ with dimensions strictly smaller than $\T X$ { which results in a low rank factorization}. Thus, we can learn a Tucker decomposition where the matrices $\m V^{(\ell)}$ are not square. Note that, in this case, the matrices define a pseudo-metric that satisfies the triangle inequality but may or may not satisfy coincidence or symmetry, depending on the $\m{\hat L}^{(\ell)}$ matrices being positive semi-definite.


	\subsection{Extension for Coupled Tensor--Matrix Recovery}\label{sec:coupled}

Many real-world tensors come with side information in the form of matrices coupled to one or more tensor modes. The above formulation of TR-MLC can incorporate such side information in the form of similarity matrices, but often the coupled matrix is not symmetric or even square. We can, however, extend the method in a relatively simple way to leverage general side information by enforcing coupled matrix-tensor decomposition. Suppose we have a tensor $\mathbf{\mc{X}}\in \mathbb R^{N_1\times \cdots\times N_c\times  \cdots\times N_K}$ and a matrix $\mathbf{M}\in \mathbb R^{J\times N_c}$ that share a common dimension $N_c$ on the $c$th tensor mode. We can extend the above framework to jointly perform Tucker decomposition of $\mc X$ and matrix factorization of $\mathbf M$ as follows: 
\begin{equation}
\begin{cases}
\mc X= \mc G\times_1 \m{V}^{(1)}\times_2 \cdots \times_{c} \m V^{(c)} \times_{c+1} \cdots \times_K \m{V}^{(K)}, \\
\m M= \m{U}^{(c)} \m{V}^{(c)},
\end{cases}
\end{equation}
where {both} Tucker and matrix decomposition contain an identical factor matrix $\m{V}^{(c)}\in \mathbb R^{n_c\times N_c}$ and $\m U^{(c)}\in \mathbb R^{J\times n_c}$ is an additional factor specific to the coupled matrix. The shared factor matrix $\m{V}^{(c)}$ ensures that the coupled matrix informs and constrains the tensor decomposition in our bilevel optimization problem. 

\section{ADMM Algorithm for TR-MLC}\label{sec:ADMM}
Bilevel optimization problems of the form~\eqref{eq:bilTen} were first introduced by Bracken and McGill \citep{bracken1973mathematical}. Later, a more general form of problem~\eqref{eq:bilTen} involving joint constraints of variables in the upper and lower problems was considered in \citep{aiyoshi1981hierarchical}. There are {many approaches} for solving general constrained bilevel optimization problems; {see, e.g.,~\citep{liu2021investigating} for survey}.  {Motivated by~\citep{hansen1992new}, we first rewrite the bilevel problem \eqref{eq:bilTen} as a single-level constrained optimization problem, and then develop} a linearized ADMM approach to solve the resulting problem. To do so, we need to calculate the derivatives of the lower-level objective given in~\eqref{eq:bilTen}. 

For brevity, we refer to the upper-level objective {in \eqref{eq:bil:up:Ten}} as $F(\T{Z})$ and the lower-level objective {in \eqref{eq:bil:low:Ten}} function as $f(\wh{\T{Z}})$ throughout this section. We also let $L(\mc{Z}; \mc{X})$ be an empirical loss function known a priori, such as the negative log-likelihood function from the generating distribution or more general objective function that depends on an unknown parameter $\mc{Z} \in \mb{R}^{N_1\times \cdots \times N_K}$~\citep{hong2020generalized,han2020optimal}. In our experiments, we used $L(\mc{Z}; \mc{X})= \frac{1}{2} \| \T X-\T Z\|^2_F$, but our derivation allows for more general losses.
\begin{lemma}[Partial Gradients of Lower-Level Objective]\label{lem:Tpgrad}
	\begin{subequations}
		\begin{eqnarray}\label{eq:Tpgrad}
		\nabla_{\hat{\T{X}}} f(\wh{\mc{Z}})= \hat{\T{M}}\times_1\hat{\m{L}}^{(1)}\times_2\cdots \times_K \hat{\m{L}}^{(K)},   \\
		\nabla_{\hat{\m{V}}^{(\ell)}} f(\wh{\mc{Z}})= \big[\hat{\m{M}}_{(\ell)} \big(\bigotimes\limits_{i\neq \ell}^K  \hat{\m{L}}^{(i)}\big) {\hat{\m{X}}_{(\ell)}^\top}    +\lambda_\ell\,\hat{\m{L}}^{(\ell)}\, \m{S}_{(\ell)}\big]\hat{\m V}^{(\ell)^\top},
		\end{eqnarray}
	\end{subequations}
	for $\ell=1, \cdots, K$. 
	Here, $\hat{\m{X}}_{(\ell)}$ denotes  $\ell$--unfolding of {$\hat{\mc{X}}$, and} 
	\begin{equation}\label{eq:mcM}
	\hat{\mc{M}}= \hat{\T{X}}\times_1 \hat{\m{L}}^{(1)}\times_2  \hat{\m{L}}^{(2)}\times_3\cdots\times_K \hat{\m{L}}^{(K)} .
	\end{equation}
\end{lemma}
\begin{proof}
	See Appendix, Section~\ref{app:lemma:proof}.
\end{proof}
Using Lemma~\ref{lem:Tpgrad}, we next show that
the proposed Problem~\eqref{eq:bilTen} can be {rewritten} as a {single-level} nonlinear constrained optimization problem. More specifically, using the optimality conditions of the lower-level Problem \eqref{eq:bil:low:Ten}, we obtain
\begin{subequations}\label{eq:bilTen:sys}
	\begin{align}
	\label{eq:sing:obj}
	\min\limits_{\mc{Z} }&\quad  F(\mc{Z}) = L( \mc{Z}; \mc{X})  +J_1(\T{G})+ \sum\limits_{\ell=1}^K J_{(2,\ell)} ({\m{V}}^{(\ell)}) \\
	\label{eq:sing:con1}
	\textrm{s.t.} & \qquad \T{M}\times_1\m{L}^{(1)}\times_2\cdots \times_K \m{L}^{(K)}=0,  \\
	& \qquad  
	\label{eq:sing:con2}
	\big[\m{M}_{(\ell)} \big(\bigotimes\limits_{i\neq \ell}^K  \m{L}^{(i)}\big) \m{X}_{(\ell)}^\top
	+\lambda_\ell\,\m{L}^{(\ell)}\, \m{S}_{(\ell)}\big]{\m V}^{(\ell)^\top}=0,  
	\end{align}
\end{subequations}
where $\T M$ is given in \eqref{eq:mcM} and $\m L^{(\ell)}=\m V^{(\ell)^\top} \m V^{(\ell)}$ for $\ell=1, \cdots, K$. One will notice that \eqref{eq:sing:con1} and \eqref{eq:sing:con2} correspond to the first-order optimality conditions of the lower-level objective in \eqref{eq:bil:low:Ten}. 

For notational simplicity, we set $K=3$ throughout this section, though the argument generalizes readily to $K>3$. In order to solve ~\eqref{eq:bilTen:sys}, we develop a linearized ADMM algorithm. ADMM~\citep{boyd2011distributed} is an attractive approach for this problem because it allows us to decouple some of the terms in ~\eqref{eq:bilTen} that are difficult to optimize jointly. 

Let 
\begin{equation}\label{eq:abr:ab}
\begin{cases}
\T A:= \T{M}\times_1\m{L}^{(1)}\times_2\m{L}^{(2)}\times_3 \m{L}^{(3)}, \\
\T B^{(\ell)}:=\big[\m{M}_{(\ell)} \big(\bigotimes\limits_{i\neq \ell}^K  \m{L}^{(i)}\big) \m{X}_{(\ell)}^\top
+\lambda_\ell\,\m{L}^{(\ell)}\, \m{S}_{(\ell)}\big]{\m V}^{(\ell)^\top}
\end{cases}
\end{equation}
and define $\T Y:= \left\{\T{Y}^{(1)}, \{\T{Y}^{(2,\ell)}\}_{\ell=1}^3 \right\}$ to be the set of dual variables corresponding to the constraints~ \eqref{eq:sing:con1} and \eqref{eq:sing:con2}.

Under this setting, the scaled augmented Lagrangian~\citep{boyd2011distributed} for 
~\eqref{eq:bilTen:sys} takes the form
\begin{multline}\label{eq:aug:T}
\mc{L} (\Upsilon)=
L(\T Z; \T X)+J_1(\T{G}) 
+ \sum\limits_{\ell=1}^3 J_{(2,\ell)}\, (\m{V}^{(\ell)})
-\Big(\langle \, \T{Y}^{(1)}\, , \, \T A \rangle\\+ \sum\limits_{\ell=1}^3 \langle \, \T{Y}^{(2,\ell)}\,, \, \T B^{(\ell)} \rangle \Big)
+ \frac{\rho}{2}\Big( \| \T A \|_F^2
+\sum\limits_{\ell=1}^3\, \|\T B^{(\ell)} \|_F^2, \Big),
\end{multline}
where 
\begin{equation}\label{eq:set}
\Upsilon:= \left\{\m{V}^{(1)}, \m V^{(2)}, \m V^{(3)}, {\T{G}}, \T Z,\,\T{Y} \right\}.
\end{equation}

Standard ADMM \citep{boyd2011distributed} applied to~\eqref{eq:aug:T} requires solving for each variable while the others are fixed; see, ~\eqref{eq:admm:upv1}--\eqref{eq:admm:upYl:T} in the Appendix, Section~\ref{app:admm}. It can be easily seen that the constrained problem~\eqref{eq:bilTen:sys} is non-convex; hence, the global convergence of standard ADMM is a priori not guaranteed. Besides, the standard ADMM requires an exact solution of each subproblem which is prohibitive in large-scale tensor completion tasks. To overcome these issues, we next develop a linearized ADMM algorithm for ~\eqref{eq:bilTen} with guaranteed convergence.
To do so, we rewrite \eqref{eq:aug:T} as
\begin{equation}\label{eq:aug:equi:T1}
\T{L} (\Upsilon) =
L(\T Z; \T X)+J_1(\T{G})+\sum\limits_{\ell=1}^3 J_{(2,\ell)} (\m{V}^{(\ell)})+ \bar{\T{L}} (\Upsilon),
\end{equation}
where 
\begin{equation}\label{eq:quadprob:T2}
\bar{\mc{L}} (\Upsilon):= \frac{\rho}{2}\Big(\| \T A-\frac{1}{\rho}\, \T Y^{(1)} \|_F^2+ \sum\limits_{\ell=1}^3\, \|\T B^{(\ell)}-\frac1\rho\, \T Y^{(2,\ell)} \|_F^2\Big).
\end{equation}
Now, using~\eqref{eq:aug:equi:T1}, we approximate each subproblem by linearizing the smooth terms in~\eqref{eq:aug:equi:T1} with respect to the factor matrices and core tensor; see, Appendix, ~\eqref{eq:reg3.3:T1}-- ~\eqref{eq:reg3.3:T4} for further details.
\begin{algorithm}[t]
	\caption{TR-MLC via Linearized ADMM} 
	\label{alg:admm1}
	\begin{algorithmic}[1]
		\State {\textbf{Input:} $\T{X} \in \mathbb{R}^{N_1 \times \cdots \times N_K}$; positive constants $\{ \lambda_{\ell}\}_{\ell=1}^K$; factor matrices $\{{\m V}^{(\ell)}_0 \in \mathbb{R}^{n_\ell\times N_\ell}\}_{\ell=1}^K$; 
			dual variables $\T Y_0$; similarity matrices $\{\m S_{(\ell)}\}_{\ell=1}^K$; number of iterations $iter$; and proximal parameters $\rho$, $\varrho^g$, $\{\varrho^\ell\}_{\ell=1}^K$}. 
		\State \textbf{Initialize:}  $t=0$,  $
		{\m V}^{(\ell)}_t ={\m V}^{(\ell)}_0$,
		$\m L^{(\ell)}_t ={\m V}^{(\ell)^\top}_t\, {\m V}^{(\ell)}_t\,$,
		$\T{Z}_t = \T{X}$, 
		
		$\wh{\T{X}}_t=\T{X} \bigotimes\limits_{\ell=1}^K {\m L}^{(\ell)}_t$, $\T G_t=\T X \bigotimes\limits_{\ell=1}^K \m V_t^{(\ell)^\top}$, and $\T{Y}_t= \T{Y}_0$. 
		\\
		\textbf{For} $t=1,2, \ldots, iter \textbf{ do}$
		\begin{itemize}
			
			\item For $\ell=1,\cdots, K$, {set} 
			$\m L^{(\ell)}:= \m V^{(\ell)^\top} \m V^{(\ell)}\in \mathbb R^{N_\ell\times N_\ell},$\\
			Update $\wh{\T X}:$
			\begin{itemize}
				\item Update factor matrices  $\m L^{(\ell)}_{t+1}$ using side/similarity matrices $\{ \m S_{(\ell)}\}_{\ell=1}^K$ and distance metric learning (following Algorithm~\ref{alg:ccA}, Appendix, Section~\ref{app:kn}) 
			\end{itemize}
			\item Update $\m V^{(\ell)}$: 
			\begin{itemize}
				\item $
				\m V^{(\ell)}_{t+1} = \text{prox}^{J_{(2,\ell)}}_{\varrho^\ell} 
				\left(\m V^{(\ell)}_{t} -\frac{1}{\varrho^\ell}  \nabla_{\m V^{(\ell)}} {\T L}(\Upsilon), \frac{1}{\varrho^\ell} \right).$
			\end{itemize}
			
			\item  Update $\T G $:
			\begin{itemize}
				\item $\T G_{t+1}=\text{prox}_{\rho_g}^{J_{1}}
				\left( \T G_t - \frac{1}{\varrho^g}\, \nabla_{\T  G} {\T L}(\Upsilon), \frac{1}{\varrho^g}\right).$
			\end{itemize}
			\item  Update $\T Z$:
			\begin{itemize}
				\item $\T Z_{t+1}\leftarrow \T G_{t+1}\times_1 \m V^{(1)}_{t+1}\times_1\m V^{(2)}_{t+1} \cdots \times_K \m V^{(K)}_{t+1}.$
			\end{itemize}
			\item Update $\T Y$:
			\begin{itemize}
				\item  $\T{Y}^{(1)}_{t+1} = \T{Y}^{(1)}_{t} - \rho\T A^{(1)}_{t+1}. $
				\item For $\ell=1,\cdots,K$: $\T{Y}^{(2,\ell)}_{t+1} = 
				\T{Y}^{(2,\ell)}_{t} - \rho \T B^{(\ell)}_{t+1}$.
			\end{itemize}
		\end{itemize}
		\textbf{End}
		\\
		\textbf{Return} $\T Z$, $\T Y$, $\T G$, and $\{\m V^{(\ell)}\}_{\ell=1}^K$.
	\end{algorithmic}
\end{algorithm}

A schematic description of the proposed linearized ADMM is given in Algorithm~\ref{alg:admm1}. The algorithm alternatively updates each variable using the proximal gradient method; see, Appendix, Section~\ref{app:admm}. In each step, it uses the proximal operator which is formally defined as: 
\begin{equation*}
\text{prox}_{t}^{J}\left(\m u\right) := \argmin \left\{ J\left(\m v\right) + \frac{t}{2}\|{\m v - \m u}\|^{2} :
\; \m v \in \R^{d} \right\},~~ \left(t > 0\right).
\end{equation*}
Here,  $J: \R^{d} \rightarrow \left(-\infty , \infty\right]$ is a proper lower semicontinuous function and $t$ is some regularization parameter. 

Note that the ADMM updates for the coupled case are very similar, requiring only slight changes to the linearized ADMM updates. These changes ensure that $\m{V}^{(c)}$ and $\m U^{(c)}$ are updated jointly. 

\section{Convergence Analysis}\label{sec:globconv}
In this section, we establish the global convergence of Algorithm~\ref{alg:admm1}. To do so, we make the following assumptions:
\begin{enumerate}[label={\textbf{(A\arabic*})}]
	\item \label{assu:pen} 
	The regularization terms $J_{1} : \R^{R_{1} \times R_{2} \cdots  \times R_{K} } \rightarrow \left(-\infty , \infty\right]$ and
	$J_{(2,\ell)} : \R^{N_{\ell} \times R_{\ell}} \rightarrow \left(-\infty , \infty\right]$ 
	are proper and lower semi-continuous such that $\inf_{\R^{R_{1} \times R_{2} \cdots  \times R_{K} }}J_{1} > -\infty$, and
	$\inf_{\R^{N_{\ell} \times R_{\ell}}} J_{(2,\ell)} > -\infty$ for $\ell=1,2 , \cdots, K$. 
	\item The loss function \label{assu:loss} 
	$L(\T{Z}; \T{X} ): \mathbb R^{N_{1} \times N_{2} \times \cdots \times N_{K}} \rightarrow \mathbb R$ is differentiable and  $\inf_{\mathbb R^{N_{1} \times N_{2} \times, \cdots, \times N_{K}}} L > -\infty$. 
	\item \label{assu:loss:lip} 
	The gradient $\nabla L(\T{Z}; \T{X})$ is Lipschitz continuous with moduli $\Pi_L$, i.e.,
	$\|\nabla L(\T{Z}_1; \T{X} ) - \nabla  L(\T{Z}_2; \T{X} )\|_F^2 \leq \Pi_L\| \T{Z}_1 - \T{Z}_2\|_F^2, $ 
	for all $\T{Z}_1, \T{Z}_2.$
\end{enumerate}
{These assumptions are common in the tensor factorization and composite optimization literature \citep{rockafellar2009variational,bauschke2011convex,tarzanagh2019regularized}.}
Note that it is not very restrictive to require  $\left\{J_1,\{ J_{(2,\ell)}\}_{\ell=1}^K\right\}$ to be proper and lower semi-continuous functions. In fact, many penalty functions including the $\ell_1$-norm, $\ell_2$-norm, $\ell_{\infty}$-norm, and nuclear norm satisfy Assumption~\ref{assu:pen}. Further, {Assumptions~\ref{assu:loss} and \ref{assu:loss:lip} cover many statistically-motivated losses used in practice, including the quadratic and exponential family losses \citep{tarzanagh2019regularized}. 
}
The crux of the proof is that the augmented Lagrangian used in the ADMM algorithm satisfies the Kurdyka--~{\L}ojasiewicz property, which is sufficient to guarantee global convergence.

For any proper, lower semi-continuous function $g : H \rightarrow (-\infty, \infty]$, we let $\partial_L g : H \rightarrow 2^H$ denote the \textit{limiting subdifferential} of $g$; see~\citep{rockafellar2009variational}~[Definition 8.3]. For any $\eta \in (0, \infty)$, we let $F_\eta$ denote the class of concave continuous functions $\varphi : [0, \eta) \rightarrow \mathbb{R}_+$ for which $\varphi(0) = 0$; $\varphi$ is $C^1$ on $(0, \eta)$ and continuous at $0$; and for all $s \in (0, \eta)$, we have $\varphi'(s) > 0$.

\begin{definition}[Kurdyka--~{\L}ojasiewicz Property]\label{defn:kl} A function $g : H \rightarrow (- \infty, \infty]$ has the \textit{Kurdyka-{\L}ojasiewicz} (KL) property at $\overline{u} \in \text{dom}(\partial_L g)$ provided that there exists $\eta \in (0, \infty)$, a neighborhood $U$ of $\overline{u}$, and a function $\varphi \in F_\eta$ such that 
	\begin{align*}
	\left(\forall u \in U \cap \{ u' \mid g(\overline{u}) < g(u') < g(\overline{u}) + \eta\}\right), \\
	\varphi'(g(u) - g(\overline{u}))\, \textnormal{dist}(0, \partial_L g(u)) \geq 1.  
	\end{align*}
	The function $g$ is said to be a \textit{KL function} provided it has the KL property at each point $u \in \text{Dom}(g)$.
\end{definition}

{Next, we present global convergence of Algorithm \ref{alg:admm1}. }
\begin{theorem}
	[Global Convergence]\label{thm:main} Suppose Assumptions~\ref{assu:pen}--\ref{assu:loss:lip} hold and the augmented Lagrangian 
	$\T{L} ({\Upsilon})$ is a KL function. Then, the sequence $\Upsilon_t= (\m V^{(1)}_t, \dots, \m V^{(K)}_t, \T{G}_t, \T{Z}_t, \T{Y}_t)$ generated by Algorithm~\ref{alg:admm1} from any starting point converges to a stationary point of \eqref{eq:aug:T}. 
\end{theorem}
We first outline the following lemma which is needed in order to prove Theorem~\ref{thm:main}. We give a formal definition of the limit point set. Let the sequence $\{\Upsilon_t\}_{t\geq 0}$ be a sequence generated by the 
	Algorithm \ref{alg:admm1} from a starting point $\Upsilon_0$. The set of all limit points is denoted by $\T C(\Upsilon_0)$, i.e.,
	\begin{equation}\label{set:lim}
	\T C(\Upsilon_0)= \left\{ \bar{\Upsilon}: \exists
	\text{ an infinite sequence } \{\Upsilon_{k_s}\}_{s \geq 0}
	\text{ s.t. }\Upsilon_{k_s} \rightarrow \bar{\Upsilon} \text{ as }~s\rightarrow\infty \right\}.
	\end{equation}
	We now show that the set of accumulations points of the sequence  $\{\Upsilon_k\}_{k\geq 0}$ generated by Algorithm \ref{alg:admm1} is nonempty and it is a subset of the critical points of $\Lc$.

	 \begin{lemma}\label{lem:8}
		Let $\{\Upsilon_t\}_{t\geq 0}$ be a sequence generated by Algorithm \ref{alg:admm1}. Then,
		\\
		(i) $\T C(\Upsilon_0)$ is a non-empty set, and any point in $\T C(\Upsilon_0)$ is a critical point of $\Lc(\Upsilon)$;\\
		(ii) $\T C(\Upsilon_0)$ is a compact and connected set;  \\
		(iii) The function $\Lc(\Upsilon)$ is finite and constant on  $\T C(\Upsilon_0)$.
	\end{lemma}
	\begin{proof}
		The proof of Lemma~\ref{lem:8} is similar to that of Lemma~15 in \citep{tarzanagh2019regularized} and is included in Appendix~\ref{app:p15}. 
		\end{proof}
{Note that the global convergence of Algorithm~\ref{alg:admm1} requires $\mc{L} (\Upsilon)$ to satisfy the KL property. In the following, we show that this assumption is not restrictive and many tensor objectives are indeed KL functions. To do so, we first introduce the definitions of semi-algebraic and sub-analytic functions  \citep{bolte2014proximal,tarzanagh2019regularized}. 
}
\begin{definition}[Semi-Algebraic Functions]\label{def:semialge}
	A function $\Psi : H \rightarrow (0, \infty]$ is \emph{semi-algebraic} provided that the graph $ \mb{G}(\Psi) = \{(x, \Psi(x)) \mid  x \in H\}$ is a semi-algebraic set, which in turn means that there exists a finite number of real polynomials $g_{ij}, h_{ij} : H \times \mb{R}\rightarrow \mb{R}$ such that 
	\begin{align*}
	\mb{G}(\Psi) := \bigcup_{j=1}^p \bigcap_{i=1}^q \{ u \in H \mid g_{ij}(u) = 0 \text{ and } h_{ij}(u) < 0\}.
	\end{align*}
\end{definition}

\begin{definition}[Sub-Analytic Functions]\label{def:analytic}
	A function $\Psi : H \rightarrow (0, \infty]$ is \emph{sub-analytic} provided that the graph $ \mathbb{G}(\Psi) = \{(x, \Psi(x)) \mid  x \in H\}$ is a sub-analytic set, which in turn means that there exists a finite number of real analytic functions $g_{ij}, h_{ij} : H \times \mathbb{R}\rightarrow \mathbb{R}$ such that 
	\begin{align*}
	\mathbb{G}(\Psi) := \bigcup_{j=1}^p \bigcap_{i=1}^q \{ u \in H \mid g_{ij}(u) = 0 \text{ and } h_{ij}(u) < 0\}.
	\end{align*}
\end{definition}

It can be easily seen that both real analytic and semi-algebraic functions are sub-analytic. Even though the sum of two sub-analytic functions may not be sub-analytic, if at least one function maps bounded sets to bounded sets, then their sum is also sub-analytic \citep{bolte2014proximal}. The KL property holds for a large class of functions including sub-analytic and semi-algebraic functions such as indicator functions of semi-algebraic sets, vector (semi)-norms $\|\cdot \|_\mu$ with $\mu \geq 0$ be any rational number, and matrix (semi)-norms (e.g., operator, trace, and Frobenius norm). These function classes cover most of the smooth and nonconvex objective functions encountered in practical applications; see \citep{bolte2014proximal,tarzanagh2019regularized} for a comprehensive list.

\begin{lemma}\label{lem:12}
The regularized augmented Lagrangian function $ \T{L} (\Upsilon)$ satisfies the KL property.
\end{lemma}

\begin{proof}
		The penalty functions $\{J_1, J_{(2, \ell)}\}_{\ell=1}^K $ satisfying Assumption~\ref{assu:pen} are semi-algebraic functions, while the loss function $L$ is sub-analytic. Further, the proximal term $\bar{\mc{L}} (\Upsilon)$ in \eqref{eq:quadprob:T2} is an analytic function. Hence, the regularized augmented Lagrangian function $ \T{L} (\Upsilon)$ given in \eqref{eq:aug:equi:T1}, which is the summation of semi-algebraic functions, is itself semi-algebraic and thus
		satisfies the KL property. 
\end{proof}

\section{Experimental Results}\label{sec:ER}
To evaluate our tensor recovery method, we compare it with the following tensor decomposition and tensor completion methods: \begin{itemize}
		\item Simultaneous tensor decomposition and completion {using
	factor prior} (STDC \citep{chen2013simultaneous}) is a constrained optimization problem that exploits rank
		minimization techniques and decomposition with the
		Tucker model. As the model structure is implicitly
		included in the Tucker model, the model uses factor priors to characterize the
		underlying distribution induced by the model factors. 
		\item High accuracy low rank tensor completion (HaLRTC \citep{liu2012tensor}) is a low rank completion method which utilizes the trace norm as a convex relaxation of an otherwise non-convex optimization problem. 
		\item Tensor Completion Made Practical,  \citep{liu2020tensor} is a recent state-of-the-art algorithm for tensor completion that involves a modification to standard alternating minimization. The method mainly tries to overcome the fact that the standard alternating minimization
		algorithm does get stuck in local minima when completing tensors with correlated components. We refer to this method as TCMP. 
		\item  Coupled Tensor Matrix Factorization (CTMF)  \citep{acar2014structure} is an algorithm for tensor decomposition that utilizes side information in form of matrices that share one mode with the tensor. The side information is supposed to provide extra features that improves the performance of the decomposition task.  
	\end{itemize}
\subsection{Datasets}
\begin{itemize}
\item The first dataset we use is the CMU face database \citep{sim2003bs}, which contains face images from 65 subjects with 11 poses and 21 types of illumination. All facial images are aligned by eye position and cropped to $32\times 32$ pixels. Images are then vectorized, and the dataset is represented as a fourth-order tensor with dimensions $(65 \times 11 \times 21\times 1024)$.

\item We also consider the Cine Cardiac dataset \citep{lingala2011accelerated}. Cine sampling is gated to a patient's heart beat, and as each data measurement is captured, it is associated with a particular cardiac phase. The dataset is represented as a fourth-order tensor with dimensions ($192\times 192\times 8\times 19$), in which $n=192, t=19$ denote a bSSFP (balanced steady-state free precession, one of the imaging techniques upon which this dataset is created) to long-axis cine sequencing parameters along with $l=8$ as the number of channels through which the signals are received. For full details of how the dataset was created, we refer the reader to \citep{candes2013unbiased}. The missing rate is $85\%$. 

\item We perform additional experiments on a public Face dataset, which is a third-order tensor with dimensions $192\times 168\times 64$ \citep{Face} and 70\% missing entries. 

\item The Fluorescence dataset \citep{bro2005standard} contains measurements of emission and excitation spectra from six different fluorophores. These were chosen on the basis of closeness to the first-order Rayleigh scatter line and their overlap in both emission and excitation spectra. In total 405 samples were recorded. Each sample was left in the instrument and scanned five consecutive times. Every measurement consists of $136$ emission wavelengths $\times 19$ excitation wavelengths. This gives a fourth-order tensor with dimensions ($405 \times 136\times 19 \times 5$).

\item We also evaluate TR-MLC on a publicly available spatial transcriptomic dataset from the adult mouse brain. We obtained the spatial gene expression tensors from a two-dimensional tissue section \footnote{ {\footnotesize\url{https://www.10xgenomics.com/resources/datasets/mouse-brain-serial-section-1-sagittal-posterior-1-standard-1-1-0}}},\footnote{
	\url{https://github.com/kuanglab/FIST}}. The dataset is a tensor of size $32285 \times 78 \times 64$, denoting the number of genes, $x$-coordinates and $y$-coordinates. Each element in the tensor represents the number of RNA copies of a gene at a particular (x,y) location. We randomly sampled $500$ genes from the full dataset, giving a tensor with dimensions ($500 \times 78 \times 64$).
\item We also used simulated datasets to assess the performance of TR-MLC. To create simulated datasets with coupled information, we followed the method of \citep{acar2014structure}. This strategy creates coupled tensor and matrices from initial ground-truth factors. For instance, the input sizes $=[20~30~40~ 50]$ and modes $=\{[1 2 3],[1,4]\}$ create a three-way tensor of size $20\times 30\times 40$ which at the first mode is coupled to a matrix of the size $20 \times 50$. 
\end{itemize}
\subsection{Experiment Details}
For each experiment, we initialized the similarity side information matrices $\m S_{(\ell)}$ using the data itself. For instance, for a three-way given tensor $\T X\in \mathbb R^{N_1\times N_2\times N_3}$, we calculated similarity information as the pairwise distance of each mode of the tensor $\T X$. Thus, coupled to the tensor $\T X$ at mode $N_\ell$, for $\ell=1,2,3$, this gives a square matrix of the size $N_\ell \times N_\ell$ which influences the $ {\m V_\ell}$ as given in problem~\eqref{eq:aug:equi:T1}. We ran TR-MLC for $20$ and $50$ iterations in all experiments unless stated otherwise. We ran STDC for 50 iterations (default setting recommended by the STDC paper). We ran HaLRTC following the example code from the authors, which uses a convergence threshold to determine the number of iterations. We ran the CMTF method following the package instructions\footnote{\url{http://www.models.life.ku.dk/joda/CMTF_Toolbox}}. It uses either a threshold or a fixed number of iterations (50) to determine convergence, whichever is achieved first. We ran the TCMP authors' Python code\footnote{\url{https://github.com/cliu568/tensor_completion}} with the default setting provided by the authors. TCMP uses a fixed number of iterations (200 by default) to determine convergence. We ran the code with this default setting.

To underscore the difference between the STDC method (the most closely related previous approach) and TR-MLC, we performed an additional experiment where we augmented the STDC objective with a distance metric learning term similar to ~\eqref{eq:bil:low:Ten}. We refer to this method as STDC$^+$ in our experiments, and we find that STDC$^+$ always outperforms STDC, highlighting the importance of this additional term (see below for details).

\subsection{Tensor Completion Performance}
We used two metrics, fit and RSE, to evaluate completion accuracy. If $\T X$ represents the true values of a tensor and $\T{\tilde X}$ is the completed version of a partially observed tensor, the fit is defined as 
$1- \frac{\| \T X- \T {\tilde X}\|_F}{\| \T X\|_F}$. Thus, fit is a number between 0 and 1, where 1 represents perfect concordance. The RSE metric is the square root of the mean-squared error between predicted and true tensor entries; a lower value is better for this metric. 

We assessed completion performance by holding out entries from each of the CMU, Cine, Fluorescence, and Adult Mouse Brain datasets; completing the tensors; and comparing the predicted values of the missing entries with their known values. The results are reported in Table ~\ref{tab:compl}. TR-MLC achieves significantly higher fit and lower RSE than both STDC and HaLRTC on the CMU, Cine, Face, and Fluorescence datasets (Table ~\ref{tab:compl}).Notably, the modified STDC objective (STDC+) using a distance metric learning term similar to that in the TR-LMC objective significantly boosts the performance of STDC. TR-MLC is also slightly better than HaLRTC on the AMB dataset, though still much better than STDC. TCMP consistently outperforms STDC, STDC+, and HaLRTC, but TR-MLC still significantly improves upon TCMP.

\begin{figure}[H]
	\centering
	\includegraphics[width=2in]{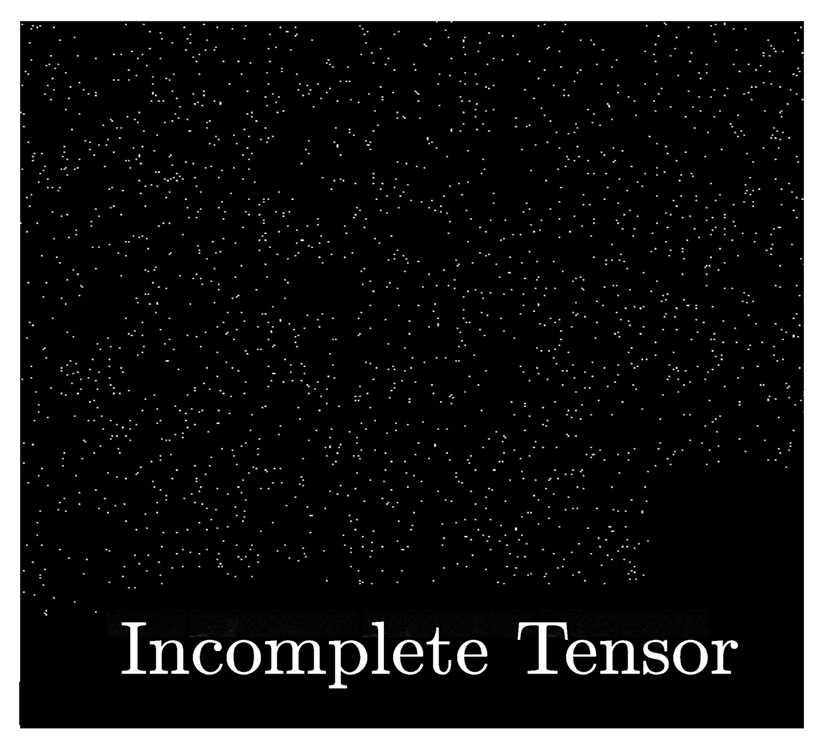}\includegraphics[width=2in]{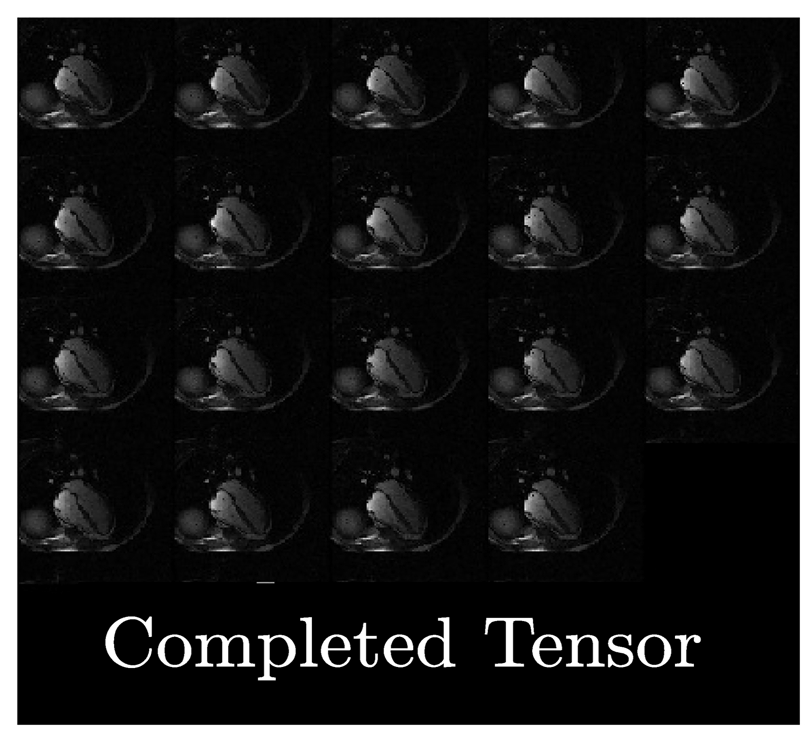}
	\caption{ The completed Cine tensor data with the initial missing rate of $85\%$ (left) using TR-MLC method after 15 iterations (right). 
	}
	\label{fig:cine}
\end{figure}

\begin{figure}
	\centering
	\includegraphics[width=3in]{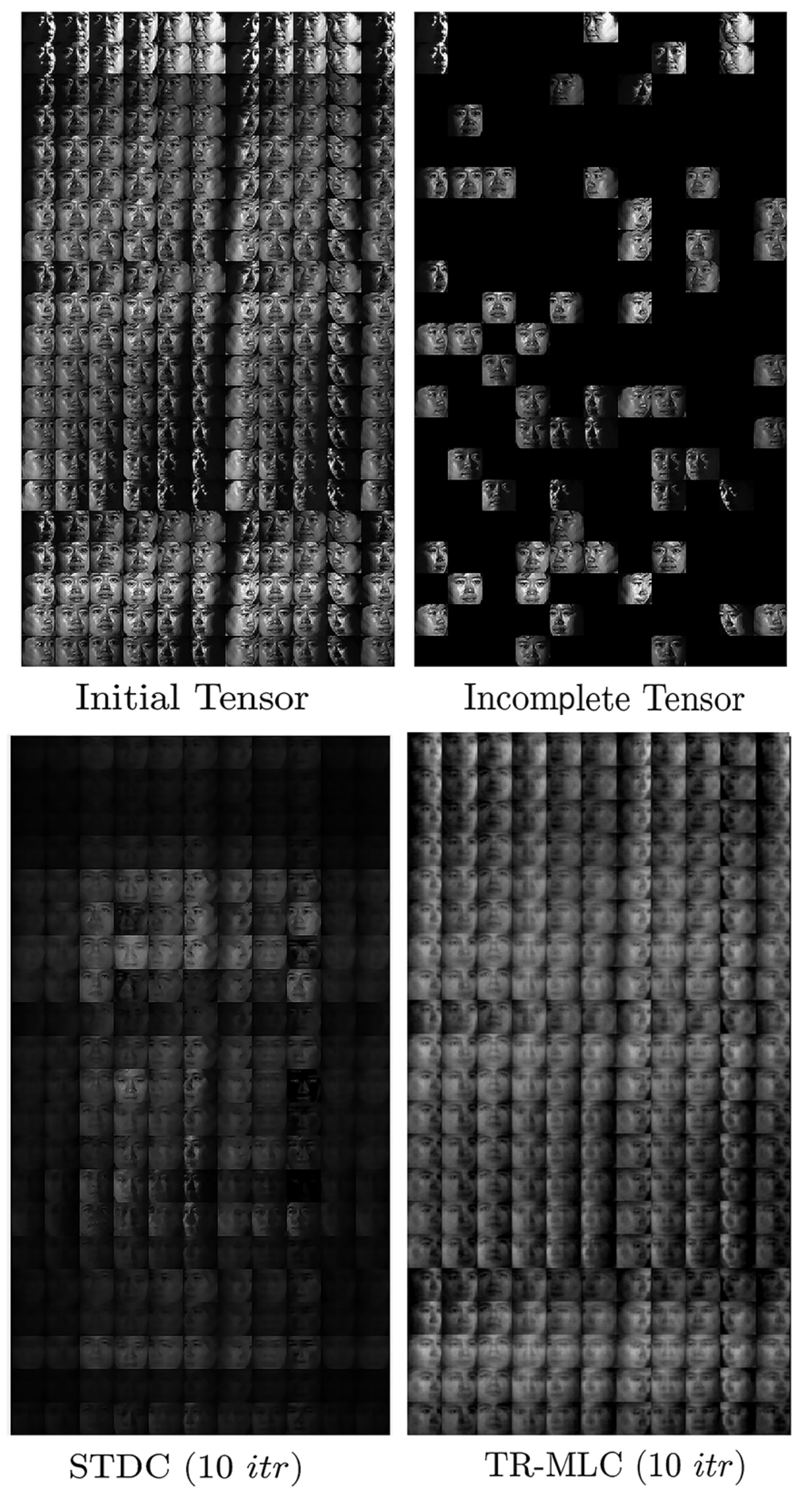}
	\caption{ Top-left shows the CMU dataset and Top-right shows the incomplete tensor with $80\%$ of the entries missing. The bottom-left and bottom-right are visualizations of the tensors completed by SDTC and TR-MLC (10 iterations), respectively. }   
	\label{fig:cmu}
\end{figure}

Qualitatively, TR-MLC is able to recover realistic-looking images of the heart from extremely sparse observations in the CINE dataset (Fig. \ref{fig:cine}). We can also see that the completed tensors (faces) from TR-MLC look qualitatively quite realistic, preserving nuances of face shape, face orientation, and lighting (Fig. \ref{fig:cmu}). In contrast, the faces reconstructed by STDC lack many of the details of the missing faces (compare left and right columns of Fig. \ref{fig:cmu}).

\begin{figure}
	\centering
	\includegraphics[width=5in]{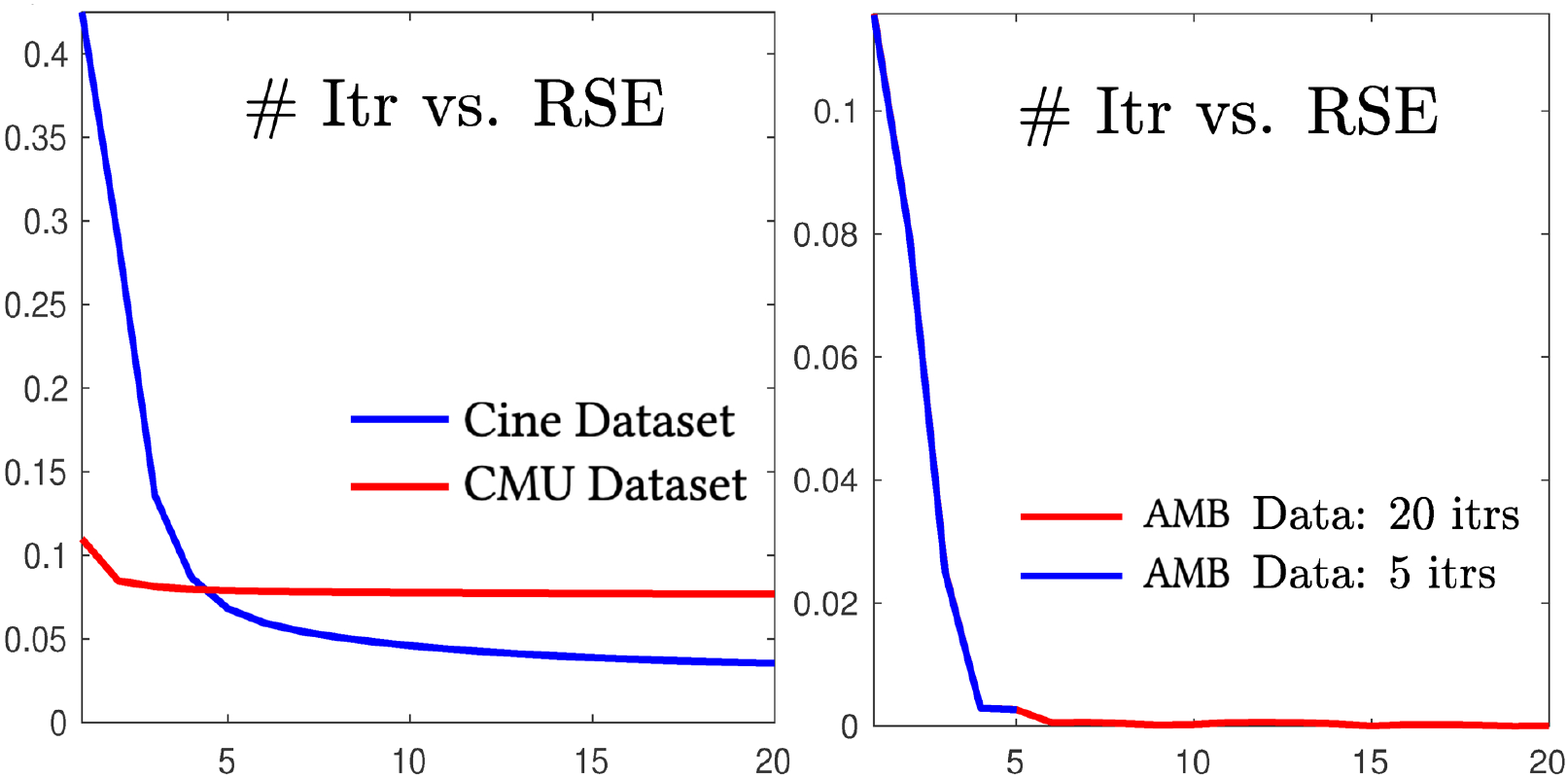}
	\caption{ (Left) Number of iterations against RSE for the 4-way data tensor Cine and CMU (see Table~\ref{tab:compl}). (Right) Same comparison for AMB data. As it shows in both figures, the RSE converges at about the 5th iteration. }
	\label{fig:rse}
\end{figure}

We also find that the ADMM algorithm for TR-MLC converges rapidly, within 5-20 iterations for the CINE, CMU, and AMB datasets (Fig. \ref{fig:rse}).

\subsection{Tensor Decomposition Performance}\label{sec:tdp}

Because our method performs simultaneous completion and decomposition, we further evaluated the decomposition performance. Unlike completion, which is simple to evaluate by holding out observed tensor entries, decomposition requires ground-truth factors, which are less often known. For this reason, we created coupled tensor and matrix data from initial factors. To this end, we a simulation approach described in \citep{acar2014structure},  and compared our results with those of CMTF (see Appendix~\ref{app:c} for further details).
Using TR-MLC, we jointly performed Tucker decomposition of the tensor, say $\T X$, and matrix decomposition of the matrix, say $\m M$, with the matrix and tensor decomposition sharing a common factor matrix along the coupled mode.
We repeated the experiment 10 times for each coupled tensor and matrix and reported the average reconstruction errors in Table \ref{tab:sim1}. 

\begin{figure}
	\centering
	\includegraphics[width=4in]{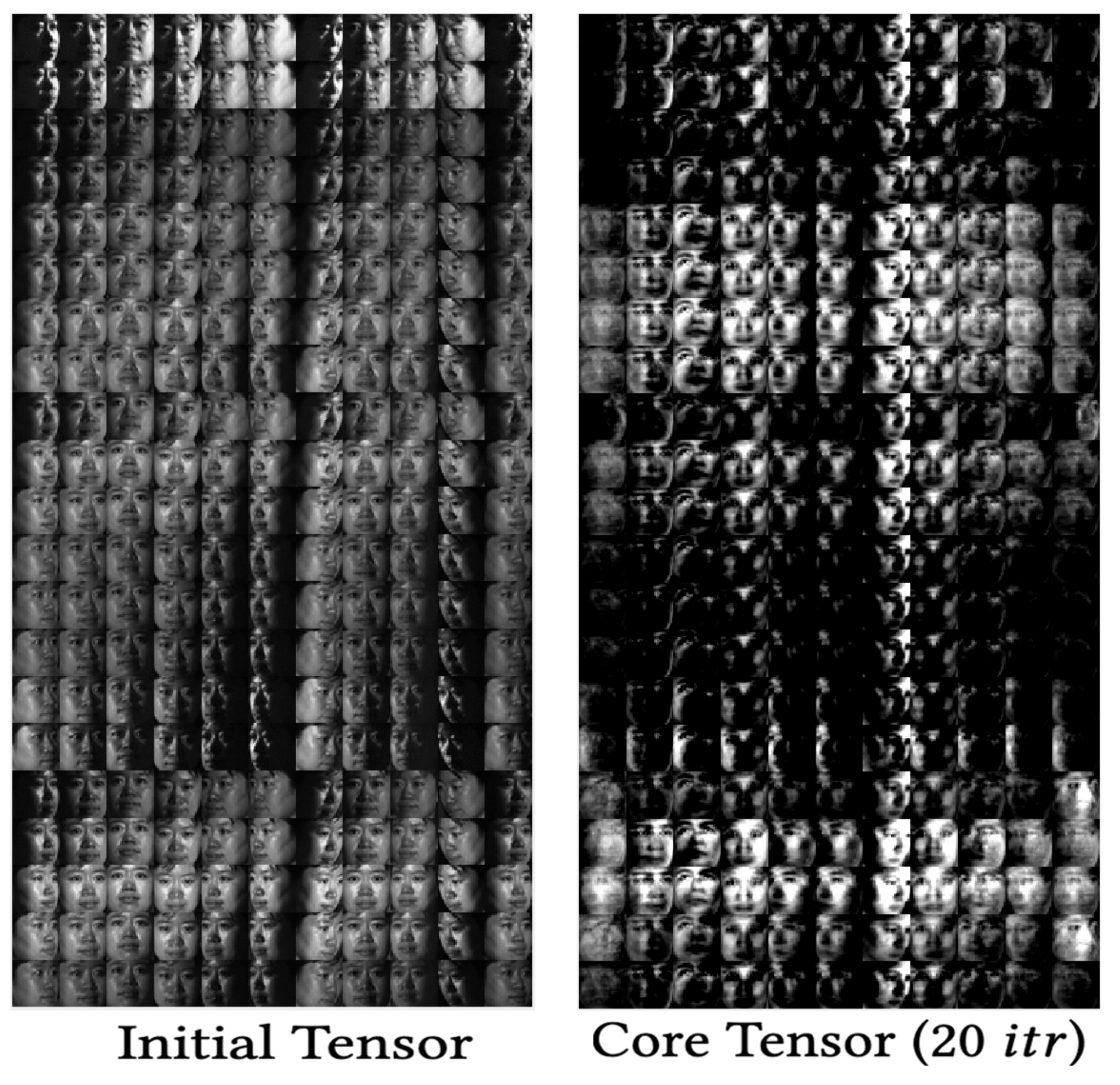}
	\caption{ This figure shows the original CMU dataset (left) and the smaller core tensor $\T G$ of size $65\times 11\times 21\times 800$ recovered by TR-MLC (right). Here, the $\# itr=20$, $fit= 0.9010$ and RSE= $0.1094$.   }
	\label{fig:cmusmall}
\end{figure}

We further utilized TR-MLC to jointly complete and decompose the CMU dataset using a core tensor whose dimensions are strictly smaller than the original tensor (Fig. \ref{fig:cmusmall}). For visualization purposes, we compressed only the pixel mode of the tensor, effectively decreasing the image resolution. Qualitatively, we see that the core tensor gives accurate but smaller representations of the original images (Fig. \ref{fig:cmusmall}). This compression comes at the cost of only a small loss in completion accuracy (fit and RSE of 0.9010 and 0.1094 with smaller core vs. fit and RSE of 0.9291 and 0.07682 with full core). 

\begin{table}
	\begin{center}
	\begin{tabular}{cccc}
		\hline \hline     
		\multicolumn{4}{c}{\textbf{Cine ($192\times 192\times 8\times 19$), 85\% missing}} \\
		\hline \hline   
		& \textbf{$\# itr$}
		&\textbf{Fit} 
		& \textbf{RSE} \\
		{\textbf{STDC}}
		& 50
		& { 0.1355}
		&{ 0.8530}\\
		{\textbf{STDC$^+$}}
		& 50
		& 0.3661
		&0.6401\\ 
		{\textbf{HaLRTC}}
		& 21
		&  0.6361
		&0.3639\\ 
		{\textbf{	TCMP}}
		& 	{{{\color{black}200 }}}
		&	{{{\color{black} 0.8121}}}
		&	{{{\color{black} 0.1903}}}\\ 
		{\textbf{TR-MLC}}
		&20
		& \m{0.9720}
		& \m{0.0354}
		\\\hline \hline
		\multicolumn{4}{c}{\textbf{CMU ($65\times 11\times 21\times 1024$), 90\% missing}}       \\ \hline \hline
		& \textbf{$\# itr$}
		&\textbf{Fit}  
		& \textbf{RSE} \\
		{\textbf{STDC}}
		&50
		&0.3837
		&0.6262\\
		{\textbf{STDC$^+$}}
		&50
		&  0.8417
		&0.1583\\  
		{\textbf{HaLRTC}}
		& 11
		&0.5829
		&0.4371\\
		{\textbf{{\color{black} 		TCMP}}}
		& 	{{{\color{black} 200}}}
		&	{{{\color{black} 0.8734}}}
		&	{{{\color{black}0.1274 }}}\\ 
		{\textbf{TR-MLC}}
		&20
		& \m{0.9291}
		&\m{0.07682 }
		\\\hline \hline
		\multicolumn{4}{c}{\textbf{Face ($192\times 168\times 64$), 70\% missing}}
		\\\hline \hline
		& \textbf{$\# itr$}
		&\textbf{Fit}  
		& \textbf{RSE} \\
		{\textbf{STDC}}
		&50
		& 0.7009
		&0.2989\\ 
		{\textbf{STDC$^+$}}
		&50
		& 0.8091
		&0.18712\\
		{\textbf{HaLRTC}}
		&43
		&  0.8162
		&0.1838\\ 
		{\textbf{{\color{black} 		TCMP}}}
		& 	{{{\color{black} 200}}}
		&	{{{\color{black} 0.8661}}}
		&	{{{\color{black} 0.1399}}}\\ 
		{\textbf{TR-MLC}}
		& 20
		&\m{0.9810}
		& \m{0.01763}
		\\\hline \hline
		\multicolumn{4}{c}{\textbf{Fluorescent ($405\times 136\times 19\times 5$), 80\% missing}} 
		\\\hline \hline
		& \textbf{$\# itr$}
		&\textbf{Fit}   
		& \textbf{RSE}                \\ 
		{\textbf{STDC}}
		&50
		&0.0002
		&1.0000\\ 
		{\textbf{STDC$^+$}}
		& 50
		&\m{0.9874}
		&\m{0.0126}\\ 
		{\textbf{HaLRTC}}
		&15 
		&0.9814
		& 0.0196\\ 
		{\textbf{{\color{black} 	TCMP}}}
		& 	{{{\color{black}200 }}}
		&	{{{\color{black} 0.9719}}}
		&	{{{\color{black} 0.0919}}}\\ 
		{\textbf{TR-MLC}}
		& 20
		&0.9833
		&0.01547\\
		\hline \hline
		\multicolumn{4}{c}{\textbf{AMB ($500\times 78\times 64$), 40\% missing}}   \\
		\hline\hline
		& \textbf{$\# itr$}
		&\textbf{Fit}   
		& \textbf{RSE} \\       
		\textbf{STDC}
		& {50}& 0.4579 & 0.7093 \\
		\textbf{STDC$^+$}
		& {50}&0.5941  & 0.3952 \\
		\textbf{HaLRTC}
		& {29}& 0.6989& 0.3152 \\
		{\textbf{{\color{black} 		TCMP}}}
		& 	{{{\color{black} 200}}}
		&	{{{\color{black}0.7035 }}}
		&	{{{\color{black} 0.2951}}}\\ 
		\textbf{TR-MLC} 
		& {20} & \m{0.7146} &\m{ 0.2666} \\
		\hline\hline
	\end{tabular}
	\end{center}
	\caption{ The method STDC$^+$ refers to the STDC method \citep{chen2013simultaneous}, where we modified it by adding the DML term to the objective function.} 
	\label{tab:compl}
\end{table}

\begin{table}[htbp]
	\begin{center}
		\begin{tabular}{c c c c   }
			\hline \hline
			\multicolumn{4}{c}{\textbf{Simulated Data}}       \\\hline \hline
			\textbf{}& \textbf{Size}  &\textbf{\# itr}  & \textbf{Avg. Err.} \\
			\textbf{CMTF} &[50 50 50 50]    & 50 &0.00695\\
			\textbf{TR-MLC} & [50 50 50 50]   & {50} &0.00606  \\
			\hline \hline
		\textbf{CMTF} & [50 50 50 50 50]    & 50 &0.08391 \\
			\textbf{TR-MLC} &[50 50 50 50 50]    &{50} &0.0691\\
			\hline \hline
				\textbf{CMTF} & [75 75 75 75] & 50 &0.08129 \\
				\textbf{TR-MLC} & [75 75 75 75]     &{50} &0.07123 \\
			\hline \hline
			\textbf{CMTF} & [100 100 100 100 100]   &50 &0.04396 \\
			\textbf{TR-MLC} &[100 100 100 100 100]  & {50} &0.04504 \\
			\hline \hline
		\end{tabular}
		\caption{ \label{tab:sim1} Here, the first three numbers in the size, e.g. $[50 50 50 50]$ denotes the size of the tensor and the numbers afterward are the size of the matrix coupled with the tensor.  The Avg. Err. denotes the average of the errors over 10 different experiments. As mentioned earlier, the first three numbers in the ``size'' indicate the size of the three-way tensor and the rest is the size of the coupled matrix. For instance, $[50~50~50~50]$ indicates a three-way tensor of the size $50\times50\times50$ to which a matrix of the size $50\times50$ is attached. The mode through which the matrix is attached is determined by another input (also see Appendix~\ref{app:c}.) }
	\end{center}
\end{table}

\section{Conclusion}\label{Conclusion}
In this manuscript, we developed a method for tensor recovery under metric learning constraints. Our approach captures the hidden features of the data by using external auxiliary information when available. It then uses a bilevel optimization scheme and ADMM to jointly complete and decompose the tensor. We have compared our results with four other methods on multiple real datasets. Moreover, our approach can also be used for simultaneous factorization when a tensor shares one (or more) modes with matrices. The bilevel optimization formulation we developed here provides a flexible framework for future extensions. In the future, we plan to use the framework to develop an Online Tensor Recovery (OTR) method for incoming time-dependent batches of information (either coupled or independent). We will also consider a variety of applications employing the proposed method, such as decision-making systems, image processing, and bioinformatics. 
\bibliography{sample}

\begin{thebibliography}{66}
\providecommand{\natexlab}[1]{#1}
\providecommand{\url}[1]{\texttt{#1}}
\expandafter\ifx\csname urlstyle\endcsname\relax
  \providecommand{\doi}[1]{doi: #1}\else
  \providecommand{\doi}{doi: \begingroup \urlstyle{rm}\Url}\fi

\bibitem[Fac()]{Face}
Publicdata.
\newblock \url{http://http://www.models.kvl.dk/datasets}.

\bibitem[Acar et~al.(2014)Acar, Papalexakis, G{\"u}rdeniz, Rasmussen, Lawaetz,
  Nilsson, and Bro]{acar2014structure}
Evrim Acar, Evangelos~E Papalexakis, G{\"o}zde G{\"u}rdeniz, Morten~A
  Rasmussen, Anders~J Lawaetz, Mathias Nilsson, and Rasmus Bro.
\newblock Structure-revealing data fusion.
\newblock \emph{BMC bioinformatics}, 15\penalty0 (1):\penalty0 1--17, 2014.

\bibitem[Aiyoshi and Shimizu(1981)]{aiyoshi1981hierarchical}
Eitaro Aiyoshi and Kiyotaka Shimizu.
\newblock Hierarchical decentralized systems and its new solution by a barrier
  method.
\newblock \emph{IEEE Transactions on Systems, Man and Cybernetics}, \penalty0
  (6):\penalty0 444--449, 1981.

\bibitem[Bagherian et~al.(2021)Bagherian, Kim, Jiang, Sartor, Derksen, and
  Najarian]{bagherian2021coupled}
Maryam Bagherian, Renaid~B Kim, Cheng Jiang, Maureen~A Sartor, Harm Derksen,
  and Kayvan Najarian.
\newblock Coupled matrix--matrix and coupled tensor--matrix completion methods
  for predicting drug--target interactions.
\newblock \emph{Briefings in bioinformatics}, 22\penalty0 (2):\penalty0
  2161--2171, 2021.

\bibitem[Bahri et~al.(2018)Bahri, Panagakis, and Zafeiriou]{bahri2018robust}
Mehdi Bahri, Yannis Panagakis, and Stefanos~P Zafeiriou.
\newblock Robust kronecker component analysis.
\newblock \emph{IEEE transactions on pattern analysis and machine
  intelligence}, 2018.

\bibitem[Bard(2013)]{bard2013practical}
Jonathan~F Bard.
\newblock \emph{Practical bilevel optimization: algorithms and applications},
  volume~30.
\newblock Springer Science \& Business Media, 2013.

\bibitem[Bauschke et~al.(2011)Bauschke, Combettes, et~al.]{bauschke2011convex}
Heinz~H Bauschke, Patrick~L Combettes, et~al.
\newblock \emph{Convex analysis and monotone operator theory in Hilbert
  spaces}, volume 408.
\newblock Springer, 2011.

\bibitem[Beckmann and Smith(2005)]{beckmann2005tensorial}
Christian~F Beckmann and Stephen~M Smith.
\newblock Tensorial extensions of independent component analysis for
  multisubject fmri analysis.
\newblock \emph{Neuroimage}, 25\penalty0 (1):\penalty0 294--311, 2005.

\bibitem[Bi et~al.(2020)Bi, Tang, Yuan, Zhang, and Qu]{bi2020tensors}
Xuan Bi, Xiwei Tang, Yubai Yuan, Yanqing Zhang, and Annie Qu.
\newblock Tensors in statistics.
\newblock \emph{Annual Review of Statistics and Its Application}, 8, 2020.

\bibitem[Bolte et~al.(2014)Bolte, Sabach, and Teboulle]{bolte2014proximal}
J{\'e}r{\^o}me Bolte, Shoham Sabach, and Marc Teboulle.
\newblock Proximal alternating linearized minimization or nonconvex and
  nonsmooth problems.
\newblock \emph{Mathematical Programming}, 146\penalty0 (1-2):\penalty0
  459--494, 2014.

\bibitem[Boyd et~al.(2011)Boyd, Parikh, and Chu]{boyd2011distributed}
Stephen Boyd, Neal Parikh, and Eric Chu.
\newblock \emph{Distributed optimization and statistical learning via the
  alternating direction method of multipliers}.
\newblock Now Publishers Inc, 2011.

\bibitem[Bracken and McGill(1973)]{bracken1973mathematical}
Jerome Bracken and James~T McGill.
\newblock Mathematical programs with optimization problems in the constraints.
\newblock \emph{Operations Research}, 21\penalty0 (1):\penalty0 37--44, 1973.

\bibitem[Bro et~al.(2005)Bro, Rinnan, and Faber]{bro2005standard}
Rasmus Bro, {\AA}smund Rinnan, and Nicolaas Klaas~M Faber.
\newblock Standard error of prediction for multilinear pls: 2. practical
  implementation in fluorescence spectroscopy.
\newblock \emph{Chemometrics and Intelligent Laboratory Systems}, 75\penalty0
  (1):\penalty0 69--76, 2005.

\bibitem[Cai et~al.(2005)Cai, He, and Han]{cai2005subspace}
Deng Cai, Xiaofei He, and Jiawei Han.
\newblock Subspace learning based on tensor analysis.
\newblock Technical report, 2005.

\bibitem[Candes et~al.(2013)Candes, Sing-Long, and Trzasko]{candes2013unbiased}
Emmanuel~J Candes, Carlos~A Sing-Long, and Joshua~D Trzasko.
\newblock Unbiased risk estimates for singular value thresholding and spectral
  estimators.
\newblock \emph{IEEE transactions on signal processing}, 61\penalty0
  (19):\penalty0 4643--4657, 2013.

\bibitem[Carroll and Chang(1970)]{carroll1970analysis}
J~Douglas Carroll and Jih-Jie Chang.
\newblock Analysis of individual differences in multidimensional scaling via an
  n-way generalization of “eckart-young” decomposition.
\newblock \emph{Psychometrika}, 35\penalty0 (3):\penalty0 283--319, 1970.

\bibitem[Chen et~al.(2013)Chen, Hsu, and Liao]{chen2013simultaneous}
Yi-Lei Chen, Chiou-Ting Hsu, and Hong-Yuan~Mark Liao.
\newblock Simultaneous tensor decomposition and completion using factor priors.
\newblock \emph{IEEE transactions on pattern analysis and machine
  intelligence}, 36\penalty0 (3):\penalty0 577--591, 2013.

\bibitem[Cichocki et~al.(2015)Cichocki, Mandic, De~Lathauwer, Zhou, Zhao,
  Caiafa, and Phan]{cichocki2015tensor}
Andrzej Cichocki, Danilo Mandic, Lieven De~Lathauwer, Guoxu Zhou, Qibin Zhao,
  Cesar Caiafa, and Huy~Anh Phan.
\newblock Tensor decompositions for signal processing applications: From
  two-way to multiway component analysis.
\newblock \emph{IEEE Signal Processing Magazine}, 32\penalty0 (2):\penalty0
  145--163, 2015.

\bibitem[Dantzig and Thapa(2006)]{dantzig2006linear}
George~B Dantzig and Mukund~N Thapa.
\newblock \emph{Linear programming 2: theory and extensions}.
\newblock Springer Science \& Business Media, 2006.

\bibitem[Davenport and Romberg(2016)]{davenport2016overview}
Mark~A Davenport and Justin Romberg.
\newblock An overview of low-rank matrix recovery from incomplete observations.
\newblock \emph{IEEE Journal of Selected Topics in Signal Processing},
  10\penalty0 (4):\penalty0 608--622, 2016.

\bibitem[De~Lathauwer et~al.(2000)De~Lathauwer, De~Moor, and
  Vandewalle]{de2000multilinear}
Lieven De~Lathauwer, Bart De~Moor, and Joos Vandewalle.
\newblock A multilinear singular value decomposition.
\newblock \emph{SIAM journal on Matrix Analysis and Applications}, 21\penalty0
  (4):\penalty0 1253--1278, 2000.

\bibitem[Dempe and Gadhi(2014)]{dempe2014new}
S~Dempe and N~Gadhi.
\newblock A new equivalent single-level problem for bilevel problems.
\newblock \emph{Optimization}, 63\penalty0 (5):\penalty0 789--798, 2014.

\bibitem[Derksen and Kemper(2015)]{derksen2015computational}
Harm Derksen and Gregor Kemper.
\newblock \emph{Computational invariant theory}.
\newblock Springer, 2015.

\bibitem[Fu and Huang(2008)]{fu2008image}
Yun Fu and Thomas~S Huang.
\newblock Image classification using correlation tensor analysis.
\newblock \emph{IEEE Transactions on Image Processing}, 17\penalty0
  (2):\penalty0 226--234, 2008.

\bibitem[Guan et~al.(2020)Guan, Dong, Absil, and Glineur]{guan2020alternating}
Yu~Guan, Shuyu Dong, P-A Absil, and Fran{\c{c}}ois Glineur.
\newblock Alternating minimization algorithms for graph regularized tensor
  completion.
\newblock \emph{arXiv preprint arXiv:2008.12876}, 2020.

\bibitem[Han et~al.(2020)Han, Willett, and Zhang]{han2020optimal}
Rungang Han, Rebecca Willett, and Anru Zhang.
\newblock An optimal statistical and computational framework for generalized
  tensor estimation.
\newblock \emph{arXiv preprint arXiv:2002.11255}, 2020.

\bibitem[Hansen et~al.(1992)Hansen, Jaumard, and Savard]{hansen1992new}
Pierre Hansen, Brigitte Jaumard, and Gilles Savard.
\newblock New branch-and-bound rules for linear bilevel programming.
\newblock \emph{SIAM Journal on scientific and Statistical Computing},
  13\penalty0 (5):\penalty0 1194--1217, 1992.

\bibitem[Hao et~al.(2019)Hao, He, Cheng, Wang, Cao, and Tao]{hao2019collect}
Fusheng Hao, Fengxiang He, Jun Cheng, Lei Wang, Jianzhong Cao, and Dacheng Tao.
\newblock Collect and select: Semantic alignment metric learning for few-shot
  learning.
\newblock In \emph{Proceedings of the IEEE/CVF International Conference on
  Computer Vision}, pages 8460--8469, 2019.

\bibitem[Hawe et~al.(2013)Hawe, Seibert, and Kleinsteuber]{hawe2013separable}
Simon Hawe, Matthias Seibert, and Martin Kleinsteuber.
\newblock Separable dictionary learning.
\newblock \emph{Computer Vision and Pattern Recognition (CVPR), 2013 IEEE
  Conference on}, pages 438--445, 2013.

\bibitem[Hillar and Lim(2013)]{hillar2013most}
Christopher~J Hillar and Lek-Heng Lim.
\newblock Most tensor problems are np-hard.
\newblock \emph{Journal of the ACM (JACM)}, 60\penalty0 (6):\penalty0 1--39,
  2013.

\bibitem[Hong et~al.(2020)Hong, Kolda, and Duersch]{hong2020generalized}
David Hong, Tamara~G Kolda, and Jed~A Duersch.
\newblock Generalized canonical polyadic tensor decomposition.
\newblock \emph{SIAM Review}, 62\penalty0 (1):\penalty0 133--163, 2020.

\bibitem[Hou and Zhang(2007)]{hou2007saliency}
Xiaodi Hou and Liqing Zhang.
\newblock Saliency detection: A spectral residual approach.
\newblock In \emph{2007 IEEE Conference on computer vision and pattern
  recognition}, pages 1--8. Ieee, 2007.

\bibitem[Kilmer et~al.(2013)Kilmer, Braman, Hao, and Hoover]{kilmer2013third}
Misha~E Kilmer, Karen Braman, Ning Hao, and Randy~C Hoover.
\newblock Third-order tensors as operators on matrices: A theoretical and
  computational framework with applications in imaging.
\newblock \emph{SIAM Journal on Matrix Analysis and Applications}, 34\penalty0
  (1):\penalty0 148--172, 2013.

\bibitem[Kolda and Bader(2009)]{kolda2009tensor}
Tamara~G Kolda and Brett~W Bader.
\newblock Tensor decompositions and applications.
\newblock \emph{SIAM review}, 51\penalty0 (3):\penalty0 455--500, 2009.

\bibitem[Kressner et~al.(2014)Kressner, Steinlechner, and
  Vandereycken]{kressner2014low}
Daniel Kressner, Michael Steinlechner, and Bart Vandereycken.
\newblock Low-rank tensor completion by riemannian optimization.
\newblock \emph{BIT Numerical Mathematics}, 54\penalty0 (2):\penalty0 447--468,
  2014.

\bibitem[Li et~al.(2021)Li, Song, Yong, and Kuang]{li2021imputation}
Zhuliu Li, Tianci Song, Jeongsik Yong, and Rui Kuang.
\newblock Imputation of spatially-resolved transcriptomes by graph-regularized
  tensor completion.
\newblock \emph{PLoS computational biology}, 17\penalty0 (4):\penalty0
  e1008218, 2021.

\bibitem[Lingala et~al.(2011)Lingala, Hu, DiBella, and
  Jacob]{lingala2011accelerated}
Sajan~Goud Lingala, Yue Hu, Edward DiBella, and Mathews Jacob.
\newblock Accelerated dynamic mri exploiting sparsity and low-rank structure:
  kt slr.
\newblock \emph{IEEE transactions on medical imaging}, 30\penalty0
  (5):\penalty0 1042--1054, 2011.

\bibitem[Liu and Moitra(2020)]{liu2020tensor}
Allen Liu and Ankur Moitra.
\newblock Tensor completion made practical.
\newblock \emph{Advances in Neural Information Processing Systems},
  33:\penalty0 18905--18916, 2020.

\bibitem[Liu et~al.(2012)Liu, Musialski, Wonka, and Ye]{liu2012tensor}
Ji~Liu, Przemyslaw Musialski, Peter Wonka, and Jieping Ye.
\newblock Tensor completion for estimating missing values in visual data.
\newblock \emph{IEEE transactions on pattern analysis and machine
  intelligence}, 35\penalty0 (1):\penalty0 208--220, 2012.

\bibitem[Liu et~al.(2021)Liu, Gao, Zhang, Meng, and Lin]{liu2021investigating}
Risheng Liu, Jiaxin Gao, Jin Zhang, Deyu Meng, and Zhouchen Lin.
\newblock Investigating bi-level optimization for learning and vision from a
  unified perspective: A survey and beyond.
\newblock \emph{IEEE Transactions on Pattern Analysis and Machine
  Intelligence}, 2021.

\bibitem[Matousek and G{\"a}rtner(2007)]{matousek2007understanding}
Jiri Matousek and Bernd G{\"a}rtner.
\newblock \emph{Understanding and using linear programming}.
\newblock Springer Science \& Business Media, 2007.

\bibitem[Milne(2017)]{milne2017algebraic}
James~S Milne.
\newblock \emph{Algebraic groups: The theory of group schemes of finite type
  over a field}, volume 170.
\newblock Cambridge University Press, 2017.

\bibitem[Narita et~al.(2012)Narita, Hayashi, Tomioka, and
  Kashima]{narita2012tensor}
Atsuhiro Narita, Kohei Hayashi, Ryota Tomioka, and Hisashi Kashima.
\newblock Tensor factorization using auxiliary information.
\newblock \emph{Data Mining and Knowledge Discovery}, 25\penalty0 (2):\penalty0
  298--324, 2012.

\bibitem[Olver and Olver(1999)]{olver1999classical}
Peter~J Olver and Peter~John Olver.
\newblock \emph{Classical invariant theory}, volume~44.
\newblock Cambridge University Press, 1999.

\bibitem[Oseledets(2011)]{oseledets2011tensor}
Ivan~V Oseledets.
\newblock Tensor-train decomposition.
\newblock \emph{SIAM Journal on Scientific Computing}, 33\penalty0
  (5):\penalty0 2295--2317, 2011.

\bibitem[Qi et~al.(2018)Qi, Shi, Sun, Wang, Yin, and Gao]{qi2018multi}
Na~Qi, Yunhui Shi, Xiaoyan Sun, Jingdong Wang, Baocai Yin, and Junbin Gao.
\newblock Multi-dimensional sparse models.
\newblock \emph{IEEE transactions on pattern analysis and machine
  intelligence}, 40\penalty0 (1):\penalty0 163--178, 2018.

\bibitem[Richardson and Slodowy(1990)]{richardson1990minimum}
Roger~Wolcott Richardson and Peter~J Slodowy.
\newblock Minimum vectors for real reductive algebraic groups.
\newblock \emph{Journal of the London Mathematical Society}, 2\penalty0
  (3):\penalty0 409--429, 1990.

\bibitem[Rockafellar and Wets(2009)]{rockafellar2009variational}
R~Tyrrell Rockafellar and Roger J-B Wets.
\newblock \emph{Variational analysis}, volume 317.
\newblock Springer Science \& Business Media, 2009.

\bibitem[Shakeri et~al.(2018)Shakeri, Bajwa, and Sarwate]{shakeri2018minimax}
Zahra Shakeri, Waheed~U Bajwa, and Anand~D Sarwate.
\newblock Minimax lower bounds on dictionary learning for tensor data.
\newblock \emph{IEEE Transactions on Information Theory}, 2018.

\bibitem[Sim(2003)]{sim2003bs}
T~Sim.
\newblock Bs, and m. bsat,“the cmu pose, illumination, and expression
  database,”.
\newblock \emph{IEEE Transactions on Pattern Analysis and Machine
  Intelligence}, 25\penalty0 (12):\penalty0 1615--1618, 2003.

\bibitem[Smilde et~al.(2005)Smilde, Bro, and Geladi]{smilde2005multi}
Age Smilde, Rasmus Bro, and Paul Geladi.
\newblock \emph{Multi-way analysis: applications in the chemical sciences}.
\newblock John Wiley \& Sons, 2005.

\bibitem[Sofuoglu and Aviyente(2020)]{sofuoglu2020graph}
Seyyid~Emre Sofuoglu and Selin Aviyente.
\newblock Graph regularized tensor train decomposition.
\newblock In \emph{ICASSP 2020-2020 IEEE International Conference on Acoustics,
  Speech and Signal Processing (ICASSP)}, pages 3912--3916. IEEE, 2020.

\bibitem[Song and Su(2021)]{song2021dstg}
Qianqian Song and Jing Su.
\newblock Dstg: deconvoluting spatial transcriptomics data through graph-based
  artificial intelligence.
\newblock \emph{Briefings in Bioinformatics}, 2021.

\bibitem[Song et~al.(2017)Song, Ge, Caverlee, and Hu]{song2017tensor}
Qingquan Song, Hancheng Ge, James Caverlee, and Xia Hu.
\newblock Tensor completion algorithms in big data analytics.
\newblock \emph{arXiv preprint arXiv:1711.10105}, 2017.

\bibitem[Su{\'a}rez-D{\'\i}az et~al.(2018)Su{\'a}rez-D{\'\i}az, Garc{\'\i}a,
  and Herrera]{suarez2018tutorial}
Juan~Luis Su{\'a}rez-D{\'\i}az, Salvador Garc{\'\i}a, and Francisco Herrera.
\newblock A tutorial on distance metric learning: Mathematical foundations,
  algorithms, experimental analysis, prospects and challenges (with appendices
  on mathematical background and detailed algorithms explanation).
\newblock \emph{arXiv preprint arXiv:1812.05944}, 2018.

\bibitem[Takeuchi and Ueda(2016)]{takeuchi2016graph}
Koh Takeuchi and Naonori Ueda.
\newblock Graph regularized non-negative tensor completion for spatio-temporal
  data analysis.
\newblock In \emph{Proceedings of the 2nd International Workshop on Smart},
  pages 1--6, 2016.

\bibitem[Tarzanagh and Michailidis(2019)]{tarzanagh2019regularized}
Davoud~Ataee Tarzanagh and George Michailidis.
\newblock Regularized and smooth double core tensor factorization for
  heterogeneous data.
\newblock \emph{arXiv preprint arXiv:1911.10454}, 2019.

\bibitem[Tucker(1964)]{tucker64extension}
L.~R. Tucker.
\newblock {T}he extension of factor analysis to three-dimensional matrices.
\newblock In H.~Gulliksen and N.~Frederiksen, editors, \emph{{C}ontributions to
  mathematical psychology.}, pages 110--127. Holt, Rinehart and Winston, New
  York, 1964.

\bibitem[Tucker et~al.(1964)]{tucker1964extension}
Ledyard~R Tucker et~al.
\newblock The extension of factor analysis to three-dimensional matrices.
\newblock \emph{Contributions to mathematical psychology}, 110119, 1964.

\bibitem[Weinberger and Saul(2009)]{weinberger2009distance}
Kilian~Q Weinberger and Lawrence~K Saul.
\newblock Distance metric learning for large margin nearest neighbor
  classification.
\newblock \emph{Journal of machine learning research}, 10\penalty0 (2), 2009.

\bibitem[Wu et~al.(2010)Wu, Zhang, and Shi]{wu2010robust}
Qiang Wu, Li-Qing Zhang, and Guang-Chuan Shi.
\newblock Robust feature extraction for speaker recognition based on
  constrained nonnegative tensor factorization.
\newblock \emph{Journal of computer science and technology}, 25\penalty0
  (4):\penalty0 783--792, 2010.

\bibitem[Xing et~al.(2002)Xing, Jordan, Russell, and Ng]{xing2002distance}
Eric Xing, Michael Jordan, Stuart~J Russell, and Andrew Ng.
\newblock Distance metric learning with application to clustering with
  side-information.
\newblock \emph{Advances in neural information processing systems}, 15, 2002.

\bibitem[Yan et~al.(2005)Yan, Xu, Yang, Zhang, Tang, and
  Zhang]{yan2005discriminant}
Shuicheng Yan, Dong Xu, Qiang Yang, Lei Zhang, Xiaoou Tang, and Hong-Jiang
  Zhang.
\newblock Discriminant analysis with tensor representation.
\newblock In \emph{2005 IEEE Computer Society Conference on Computer Vision and
  Pattern Recognition (CVPR'05)}, volume~1, pages 526--532. IEEE, 2005.

\bibitem[Yang and Jin(2006)]{yang2006distance}
Liu Yang and Rong Jin.
\newblock Distance metric learning: A comprehensive survey.
\newblock \emph{Michigan State Universiy}, 2\penalty0 (2):\penalty0 4, 2006.

\bibitem[Zhang et~al.(2014)Zhang, Ely, Aeron, Hao, and Kilmer]{zhang2014novel}
Zemin Zhang, Gregory Ely, Shuchin Aeron, Ning Hao, and Misha Kilmer.
\newblock Novel methods for multilinear data completion and de-noising based on
  tensor-svd.
\newblock In \emph{Proceedings of the IEEE Conference on Computer Vision and
  Pattern Recognition}, pages 3842--3849, 2014.

\bibitem[Zhao et~al.(2015)Zhao, Zhang, and Cichocki]{zhao2015bayesian}
Qibin Zhao, Liqing Zhang, and Andrzej Cichocki.
\newblock Bayesian cp factorization of incomplete tensors with automatic rank
  determination.
\newblock \emph{IEEE transactions on pattern analysis and machine
  intelligence}, 37\penalty0 (9):\penalty0 1751--1763, 2015.

\end{thebibliography}
\clearpage
\section{Appendix}
\subsection{Change of Coordinate Using Kempf-Ness Theorem}\label{app:kn}
Here, we use the representation theory of reductive groups \citep{milne2017algebraic,richardson1990minimum} for group $H$, a \emph{reductive} linear algebraic group. The formal definition of reductive group is given in \citep{derksen2015computational}.  
Intuitively, \emph{reductivity} can be thought of as the property of a linear algebraic group $H$ that guarantees a representation for the group $H$. 
For a real (or complex) vector space $X$, the general linear group, $\GL(X)$, and  special linear group, $\SL(X)$, are reductive groups and so are their products. In general, $\GL_d(X)$ is the set of $d\times d$ invertible matrices over $X$, together with the \textit{matrix multiplication} operation and $\SL_d(X)$ is a subset of $\GL_d(X)$ consisting of those elements whose determinants are 1. An $n\times n$ real matrix can be thought of as an element in $X\otimes X\cong \R^{n\times n},$ where $X$ and $Y$ are vector spaces of dimension $n$ and $m$, respectively. Similarly, an $n\times m\times p$ tensor is an element in the representation $X\otimes Y\otimes Z\cong \R^{n\times m\times p}$ of the reductive group  $\GL(X)\times \GL(Y)\times \GL(Z)$.
The Kempf-Ness Theorem (Theorem ~\ref{thm:KNT}) can be thought of as a non-commutative duality theory paralleling the linear programming duality given by the \emph{Farkas} lemma \citep{matousek2007understanding}, which corresponds to the commutative world.
The theorem, which is a result of the geodesic convexity of the mapping $\varphi$, implies that there is a unique metric--the Euclidean metric after the change of coordinates given by $h$.
A procedure for performing change of coordinates based on the Kempf-Ness Theorem is given in Algorithm~\ref{CCAlg}. 

\begin{algorithm}[h]
	\caption{Change of Coordinates/Kempf-Ness Normalization}
	\label{CCAlg}
	\begin{algorithmic}[1]
		\State{\textbf{Input:} Tensor $\T X$ (or a set of tensors $\T X_1, \cdots, \T X_m\in \mathbb R^{n\times p\times q}$), reg. para. $\lambda$, number of iterations $iter$}
		\State{	\textbf{Output:} New set of coordinates $(\m U_1, \m U_2, \m U_3)\in  H$ and new $\T X$}
		
		\State {\textbf{Function}~CoordinateChange($\T X_1, \T X_2, \cdots, \T X_m, \lambda, iter$)}
		\State $\m U_1\leftarrow I_{n}$, \quad $\m U_2\leftarrow I_{p}$, \quad $\m U_3\leftarrow I_{q}$
		\For{$i=1, \cdots, m$   }
		\State Concatenate all tensors and matricize: $\m X_{(n)}, \m X_{(p)},\m X_{(q)}$ 
		\EndFor
		\State minimize $\min\limits_j=\|\T X^{(j)}\|$
		\State Min$_j \leftarrow$ New $\min_j$
		\While {$\frac{\text{Min}_j}{ \min_j}<1-\lambda$ \& number of iterations is less than $iter$}
		\For {$j=1,2,3$}
		\State $\m U_j^*\leftarrow \argmin\limits_j\,\|\m U_j\cdot \T X^{(j)}\|$
		\State $\T X^{(j)}\leftarrow \m U_j^*\times_j \T X^{(j)}$
		\State min$_j\leftarrow \|\T X^{(j)}\|$
		\State $\m U_j\leftarrow \m U^*_j\cdot \m U_j$
		\EndFor
		\EndWhile\\
		\Return $(\m U_1, \m U_2, \m U_3)$ and $\T X$
		\State{\textbf{end Function}}
	\end{algorithmic}
\end{algorithm}
Next, in order to determine the optimal metric, we assume that $\m x_1,\m x_2,\dots,\m x_m$ in $X\cong \R^n$ are $m$ data points  whose respective invertible covariance matrix is shown by $\mathbf M$. Then $\mathbf x=( x_1,\dots,x_m)\in V=X^m\cong \R^{n\times m}$ and the function $\varphi:\SL(X)\to \R$, defined by $\varphi(\m U)=\|\m U\cdot \mathbf x\|^2$
has a critical point, namely $\bar {\m U}=\mathbf M^{-\frac{1}{2}}$. 
The optimal metric is indeed the Mahalanobis distance if the data points $\m x_1,\m x_2,\dots,\m x_m$ are evenly distributed. Otherwise, a better choice of $H$ yields a more optimal metric. Determining an optimal choice of $H$ induces a metric and regularization terms that are directly used 
in the algorithm. For instance, for a three-way tensor $\T X\in X\otimes Y\otimes Z$, one can optimize $
\varphi(\m U_1,\m U_2,\m U_3)=\|(\m U_1,\m U_2,\m U_3)\cdot \T X\|^2,$
for $(\m U_1,\m U_2,\m U_3)\in H:=\SL(X)\times \SL(Y)\times \SL(Z)$, using
alternating optimization. That is, we can first optimize for $\m U_1\in \SL(X)$ while fixing $\m U_2$ and $\m U_3$, followed by optimizing $\m U_2$ holding the others fixed. The process can then be iterated until the desired convergence. 
Each optimization step reduces
to the case of $m$ data points $\m x_1,\m x_2,\dots,\m x_m$ in $X\cong \R^n$ with invertible covariance matrix $\mathcal M$, which was discussed above. It can be shown that this 
procedure converges to an optimal solution and in practice only a few iterations are needed. Lastly, it is noteworthy that the choice of the group $H$ is not unique. The following algorithm shows how to perform distance metric learning using the Kemp-Ness Theorem \ref{thm:KNT}.

\begin{algorithm}[H]
	\caption{Learning Factors using DML}
	\label{alg:ccA}
	\begin{algorithmic}[1]
		\State {\textbf{Input}: A factor matrix with two modes, similarity matrices $\m S_{XX}, \m S_{YY},$ sizes of each mode: $n_X, n_Y,$ observed index set at entries  $\m\Omega_{XY}=((i_1,j_1),(i_2,j_2),\dots,(i_k,j_k)).$ Fixed vector $v_{XY}\in \mb R^k$, where $\m\Omega_{XY}(M_{XY})=\nu_{XY},$ regularization parameters at each mode $\lambda_X$ and $\lambda_Y$, and  $iter$}
		\State {\textbf{Initialize}: $\m L_X \leftarrow \mb I_{n_X}$,
			and   $\m L_Y \leftarrow \mb I_{n_Y}$}
		\For{$i=1, \cdots, iter$   }
		\State Choose $\m M_{XY}$ that minimizes $~\|\m L_X\, \m M_{XY}\, \m L_Y^\top\|^2_F~$ under constraint $\m\Omega_{XY}(\m M_{XY})=\nu_{XY}$
		\State {$\m L_X \leftarrow \big(\lambda_{X}\, \m S_{XX}\, +\m M_{XY}\, \m L_X^\top \, \m L_X\, \, \m M_{XY}^\top\big)^{-1/2}$}
		\State {$\m L_X \leftarrow \det (\m L_X)^{-1/n_X}\, \m L_X$}
		\State {$\m L_Y \leftarrow \big(\lambda_{Y}\, \m S_{YY}\, +\m M_{XY}^\top \m L_Y^\top \, \m L_Y\, \m M_{XY}\,\big)^{-1/2}$}
		\State {$\m L_Y \leftarrow \det (\m L_Y)^{-1/n_Y}\, \m L_Y$}
		\EndFor
		\State{\textbf{End}}\\
		\Return factors $\m L_X$, $\m L_Y$
	\end{algorithmic}
\end{algorithm}


\subsection{Proof of Lemma \ref{lem:pgrad:T}}\label{app:lemma:proof}

\begin{proof}
	For simplicity, let $K=3$. {From the chain rule, we have}
	\begin{equation*}
	\begin{split}
	\frac{\partial f}{\partial \hat{\m{L}}_{ij}^{(1)}} & = \sum_{a=1}^{N_1} \sum_{b=1}^{N_2} \sum_{c=1}^{N_3} \frac{\partial f}{\partial \hat{\T{Z}}_{abc}} \cdot \frac{\partial \hat{\T{Z}}_{abc}}{\partial \hat{\m{L}}_{ij}^{(1)}}  + \lambda_\ell (\hat{\m{L}}^{(\ell)}\, \m{S}_{(\ell)})_{ij} \\
	& = \sum_{a=1}^{N_1} \sum_{b=1}^{N_2} \sum_{c=1}^{N_3} \frac{\partial f}{\partial \hat{\T{Z}}_{abc}}
	\cdot \left( \mathbb I_{\{a=i\}} \sum_{r_2=1}^{R_2}\sum_{r_3=1}^{R_3} \T{Y}_{jr_2r_3}\hat{\m{L}}_{b r_2}^{(2)}\hat{\m{L}}_{cr_3}^{(3)} \right) \\
	&+ \lambda_\ell (\hat{\m{L}}^{(\ell)}\, \m{X}_{(\ell)})_{ij}\\
	& = \sum_{b=1}^{N_2} \sum_{c=1}^{N_3} \frac{\partial f}{\partial \hat{\T{Z}}_{ibc}} \left(\sum_{r_2=1}^{R_2}\sum_{r_3=1}^{R_3} \T{Y}_{jr_2r_3}\hat{\m{L}}_{br_2}^{(2)}\hat{\m{L}}_{c r_3}^{(3)}\right) 
	+ \lambda_\ell (\hat{\m{L}}^{(\ell)}\, \m{S}_{(\ell)})_{ij},
	\end{split}
	\end{equation*}
	where the second identity comes from the fact:
	$$\hat{\T{Z}}_{abc} = \sum_{r_1=1}^{R_1}\sum_{r_2=1}^{R_2}\sum_{r_3=1}^{R_3} \T{Y}_{r_1r_2r_3}\hat{\m{L}}_{ar_1}^{(1)}\hat{\m{L}}_{br_2}^{(2)}\hat{\m{L}}_{cr_3}^{(3)}.$$
	
	Let $\hat{\T{M}}= \wh{\T{X}}\times_1 \hat{\m{L}}^{(1)}\times_2  \hat{\m{L}}^{(2)}\times_3\cdots\times_K \hat{\m{L}}^{(K)}$. One can verify that
	\begin{multline*}
	\left(\hat{\T{M}}_{(1)}(\hat{\m{L}}^{(3)}
	\otimes \hat{\m{L}}^{(2)})
	\m{\wh X}_{(1)}^\top\right)_{ij} \\= \sum_{k_1=1}^{N_2}
	\sum_{k_2=1}^{N_3}\sum_{k_3=1}^{R_2} \sum_{k_4=1}^{R_3} 
	\frac{\partial f}{\partial 
		\hat{\T{Z}}_{i k_1k_2}}\cdot \hat{\m{L}}_{k_1k_3}^{(2)}\hat{\m{L}}_{ k_2k_4}^{(3)} \wh{\T X}_{j k_3k_4}.
	\end{multline*}
	Here, $\hat{\m{M}}_{(1)}$ and $\wh{\m{X}}_{(1)}$ are mode-$1$ matricization of $\hat{\T{M}}$ and $\wh{\T{X}}$, respectively.
	Now, applying \emph{chain rule} using the relation $\hat {\m L}^{(\ell)}= \hat {\m V}^{(\ell)^\top} \hat {\m V}^{(\ell)}$, the rest follows.
\end{proof}

\subsection{Derivations for ADMM Updates}\label{app:admm}
In a typical iteration of the ADMM for minimizing~\eqref{eq:aug:T}, the updates to be implemented are: 
\begin{subequations}
	\begin{eqnarray} \label{eq:admm:up:T}
	\m V^{(1)}_{t+1} &=&  
	\argmin_{\mathbf V^{(1)}} \quad  \mc{L}\big(\underbrace{
		\m V^{(1)}, \m V^{(2)}_t, \m V^{(3)}_t, {\T{G}}_t, 
		\T Z_t, \T{Y}_{t}}_{=:\Upsilon_t^{\m V^{(1)}}}\big),
	\label{eq:admm:upv1}\\
	\m V^{(2)}_{t+1}&=&  
	\argmin_{\mathbf V^{(2)}} 
	\quad \mc{L}\big(\underbrace{
		\m V^{(1)}_{t+1}, \m V^{(2)}, 
		\m V^{(3)}_t, {\T G}_t,
		\T Z_t, \T{Y}_{t}}_{=:\Upsilon_t^{\m V^{(2)}}}\big),
	\label{eq:admm:upv2}\\
	\m V^{(3)}_{t+1} &=&  
	\argmin_{\mathbf V^{(3)}} \quad  \mc{L}\big(\underbrace{
		\m V^{(1)}_{t+1}, 
		\m V^{(2)}_{t+1}, 
		\m V^{(3)},{\T G}_t,
		\T Z_t, \T{Y}_{t} }_{=:\Upsilon_t^{\m V^{(3)}}}\big),
	\label{eq:admm:upv3}\\
	{\T G}_{t+1} &=&  
	\argmin_{{\T G}} \quad  \mc{L} \big(\underbrace{\m V^{(1)}_{t+1},
		\m V^{(2)}_{t+1}, \m V^{(2)}_{t+1}, 
		{\T G},
		\T Z_t,  \T{Y}_{t}}_{=:\Upsilon^{{\T G}}_t}\big),
	\label{eq:admm:ups} \\
	{\T Z}_{t+1} &=&  
	\argmin_{{\T Z}} \quad  \mc{L} \big(\underbrace{\m V^{(1)}_{t+1},
		\m V^{(2)}_{t+1}, \m V^{(2)}_{t+1}, 
		{\T G}_{t+1},
		{\T Z},  \T{Y}_{t}}_{=:\Upsilon^{{\T Z}}_t}\big),
	\label{eq:admm:zup} \\
	\T{Y}^{(1)}_{t+1} &=&
	\T{Y}^{(1)}_{t}  -\rho\,\T A_{t+1},
	\quad 
	\T{Y}^{(2,\ell)}_{t+1}=  
	\T{Y}^{(2,\ell)}_{t}-\rho\,\T B^{(\ell)}_{t+1},
	\label{eq:admm:upYl:T}
	\end{eqnarray}
\end{subequations}
where $\T A$ and $\T B^{(\ell)}$ are given in Eq.~\eqref{eq:abr:ab}.
We regularize each subproblem in ~\eqref{eq:admm:up:T}--\eqref{eq:admm:ups} as follows:
\begin{subequations}
	\begin{eqnarray} 
	\m V^{(\ell)}_{t+1}&=&  
	\argmin_{\m V^{(\ell)}} \quad  \mc{L}(\Upsilon_t^{\m V^{(\ell)}})
	+\frac{ \varrho^\ell}{2 }\| 
	\m V^{(\ell)}- \m V^{(\ell)}_t\|_F^2,
	\label{eq:reg3.3:T1}\\
	{\T{G}}_{t+1} &=&  
	\argmin_{{\T G}} \quad  \mc{L}(\Upsilon_t^{{\T G}})
	~+ ~ \frac{ \varrho^g}{2 } \| {\T  G}-{\T G}_t\|_F^2,
	\label{eq:reg3.3:T4}
	\end{eqnarray}
\end{subequations}
where positive constants $ \varrho^g$, and $\{\varrho^\ell\}_{\ell=1}^{3}$ correspond to the  regularization parameters.  

Now, we approximate \eqref{eq:reg3.3:T1}--\eqref{eq:reg3.3:T4} by linearizing the function $\bar{\mc{L}} (\Upsilon)$ given in \eqref{eq:quadprob:T2} w.r.t. $\m V^{(1)}$, $\m V^{(2)}$, $\m V^{(3)}$, and ${\T G}$ as follows:
\begin{subequations}
	\begin{eqnarray} 
	\nonumber
	\m V^{(\ell)}_{t+1}=
	\argmin\limits_{\m V^{(\ell)}}
	\langle \nabla_{{\m V^{(\ell)}}} 
	\bar{\T{L}}(\Upsilon_t^{\m V^{(\ell)}})\, , \,
	\m V^{(\ell)}-\m V^{(\ell)}_t \rangle +\, J_{(2,\ell)}(\m V^{(\ell)})&&\\
	+\frac{\varrho^\ell}{2}
	\| \m V^{(\ell)}-\m V^{(\ell)}_t\|_F^2, && \label{eq:finalup1}\\
	{\T G}_{t+1} =       \argmin\limits_{{\T G}}
	\langle \nabla_{{{\T G}}}
	\bar{\T{L}} (\Upsilon^{{\T G}}_t)    
	\, ,\,  {\T G}-{\T G}_t \rangle +\frac{\varrho^g}{2 }\| {\T G}-{\T G}_t\|_F^2,  &&  \label{eq:finalup3}
	\end{eqnarray}
\end{subequations}
where $\nabla_{{\m V^{(\ell)}}} \bar{\mc{L}}$ and $\nabla_{{\T G}} \bar{\mc{L}}$ denote the gradients of \eqref{eq:quadprob:T2} w.r.t. $\m V^{(\ell)}$, for $\ell=1,2,3$, and ${\T G}$, respectively. 

The following two lemmas give the partial gradients of ${\Lc}(\Upsilon)$ given in \eqref{eq:quadprob:T2}
as well as the partial gradient of  $\mc{L} (\Upsilon)$ given in ~\eqref{eq:aug:equi:T1}. 
\begin{lemma} \label{lem:pgrad:T} The partial gradients of ${\Lc}(\Upsilon)$ given in \eqref{eq:quadprob:T2} w.r.t. $\wh{\T X}$ and $\m V^{(\ell)}$, for $\ell=1,2,3$ are 
	\begin{subequations}
		\begin{align}
		&\nabla_{\m V^{(\ell)}}
		\nonumber
		{\Lc}(\Upsilon) =\rho\,
		\Big[\big(\T M \times_1 \wh{\m \Lambda}\big)\big(
		\wh{\m X}_{(\ell)} 
		\wh{\m \Lambda }_{\neq \ell} \,  \wh{\m \Lambda}^\top\, \,  \wh{\m \Lambda}^\top 
		+ \wh{\m \Lambda}_{\neq \ell}\,  \T M^\top\big)\m{V}^{(\ell)^\top}\\
		\nonumber
		&+\sum\limits_{\ell=1}^3\,     \Big(\m M_{(\ell)}\, \wh{\m \Lambda}_{\neq \ell}\,  \wh{\m{X}}_{(\ell)}^\top 
		+\lambda_\ell\,\m{V}^{(\ell)^\top}\, \m{V}^{(\ell)}\, \m{S}_{(\ell)} \Big)
		\\
		&
		\Big( 
		\wh{\m X}_{(\ell)} \wh{\m \Lambda}_{\neq \ell} \, 
		\wh{\m \Lambda}_{\neq \ell}^\top\, \wh{\m S}_{(\ell)}^\top+\lambda_\ell\, \m S_{(\ell)} \Big)\m{V}^{(\ell)^\top}\Big],
		\label{eq:lemma3:1}\\
		\nonumber
		&\nabla_{\wh{\T{X}}}   
		{\Lc}(\Upsilon) = \rho\, \Big[\T M\times_1\wh{\m \Lambda} \times_1\wh{\m \Lambda}\times_1\wh{\m \Lambda}+\sum\limits_{\ell=1}^3\,     \Big(\m M_{(\ell)}\,  \wh{\m \Lambda}_{\neq \ell}\,  \wh{\m{X}}_{(\ell)}^\top \\
		&+\lambda_\ell\,
		\m{V}^{(\ell)^\top}\m{V}^{(\ell)}\, \m{S}_{(\ell)} \Big)
		\Big( \wh{\m X}_{(\ell)}\wh{\m \Lambda}^\top \, 
		\wh{\m \Lambda}_{\neq \ell}\, \wh{\m X}_{(\ell)^\top }+ \m M_{(\ell)}\, \wh{\m \Lambda}\, \Big)\m{V}^{(\ell)^\top}\Big], 
		\label{eq:lemma3:2}
		\end{align}
	\end{subequations}
	where $\T M $ given in ~\eqref{eq:mcM}, $\m M_{(\ell)}$ and $\wh{\m X}_{(\ell)}$ denote the $\ell$-th unfolding of tensor $\T M $ and $\wh{\T X}$ respectively, and 
	\begin{equation}
	\nonumber
	\wh{\m \Lambda}:= \bigotimes_{\ell=1}^3 \m L^{(\ell)},\quad ~~ \quad 
	\wh{\m \Lambda}_{\neq \ell}:= \bigotimes_{i\neq\ell}^3 \m L^{(i)}. 
	\end{equation}
\end{lemma}
{\small
	\begin{proof}
		Proof follows from Lemma.~\ref{lem:pgrad:T} and chain rule. 
\end{proof}}

The following lemma gives the partial gradient of $\mc{L} (\Upsilon)$ given in ~\eqref{eq:aug:equi:T1}
w.r.t. $\wh{\T X}$, $\T G$, and $\m V^{(\ell)}$ for $\ell=1,2,3$. 
\begin{lemma}
	The partial gradients of $\mc{L} (\Upsilon)$ w.r.t. $\wh{\T X}$, $\T G$, and $\m V^{(\ell)}$ for $\ell=1,2,3$,  are given as:
	\begin{subequations}
		\begin{align}
		\nabla_{\m V^{(\ell)}}
		\mc{L} (\Upsilon)&=\m Z_{(\ell)}\big(\bigotimes_{i\neq \ell}^3 \m V^{(i)}\big) \m G_{(\ell)}+
		\nabla_{\m V^{(\ell)}}
		{\Lc}(\Upsilon), 
		\label{eq:grad:L:1}\\
		\nabla_{\wh{\T X}}   
		\mc{L} (\Upsilon) &= \nabla_{\wh{\T X}}
		{\Lc}(\Upsilon), \quad \textnormal{and}\quad \nabla_{\T{G}}   
		\mc{L} (\Upsilon) = \T Z, 
		\label{eq:grad:L:3}
		\end{align}
	\end{subequations}
	where $\T Z=\T G\times_1 \m V^{(1)}\times_2 \m V^{(2)}\times_3 \m V^{(3)}$. Here $\m Z_{(\ell)}$ and $\m G_{(\ell)}$ are the mode-$\ell$ matricization of $\T Z$ and $\T G$, respectively.$\nabla_{\m V^{(\ell)}}
	{\Lc}(\Upsilon)$ and $\nabla_{\wh{\T X}}
	{\Lc}(\Upsilon)$ are given in Lemma.~\ref{lem:pgrad:T}.
\end{lemma}


\subsection{Proof of Lemma~\ref{lem:8}}\label{app:p15}
\begin{proof}
	(i). The sequence $\{\Upsilon_t\}_{t\geq 0}$ is bounded, which implies that $\T C(\Upsilon_0)$ is non-empty due to the Bolzano-Weierstrass Theorem. Consequently, there exists a sub-sequence $\{\Upsilon_{t_s}\}_{s\geq 0}$ such that 
	\begin{equation}\label{lim:seeq}
	\Upsilon_{t_s}\rightarrow \Upsilon_*, \qquad \text{as} \quad s\rightarrow \infty.
	\end{equation}
	Since $J_{3,n}$ is lower semi-continuous, \eqref{lim:seeq} yields
	\begin{equation} \label{eq:inf}
	\liminf_{s\rightarrow\infty} J_{2,\ell} (\m V_{t_s}^{(\ell)}) \geq  J_{2,\ell} (\m V_{*}^{(\ell)}).
	\end{equation}
	Further, from the iterative step \eqref{eq:finalup1} , we have 
	\begin{eqnarray}
	\m V^{(\ell)}_{t+1} = \argmin_{\m V^{(\ell)}} \quad \lambda_{2,\ell} J_{2,\ell}(\m V^{(\ell)})+ \langle \nabla_{\m V^{(\ell)}} \bar{\Lc}(\Upsilon_t^{\m V^{(\ell)}_t}), \m V^{(\ell)}-\m V^{(\ell)}_t \rangle + \frac{ \varrho^\ell}{2 }\| \m V^{(\ell )}-\m V^{(\ell)}_t\|_F^2.\notag
	\end{eqnarray}
	Thus, letting $\m V^{(\ell)}= \m V^{(\ell)}_{*}$ in the above, we get
	\begin{eqnarray}\label{eq:lin}
	\nonumber
	&& \lambda_{2,\ell} J_{2,\ell}(\m V^{(\ell)}_{t+1})+ \langle \nabla_{\m V^{(\ell)}} \bar{\Lc}(\Upsilon_t^{\m V^{(\ell)}_t}), \m V^{(\ell)}_{t+1}-\m V^{(\ell)}_t \rangle + \frac{ \varrho^\ell}{2 }\| \m V^{(\ell)}_{t+1}-\m V^{(\ell)}_t\|_F^2, \\
	&\leq &  \lambda_{2,\ell} J_{2,\ell}(\m V^{(\ell)}_*)+ \langle \nabla_{\m V^{(\ell)}} \bar{\Lc}(\Upsilon_t^{\m V^{(\ell)}_t}), \m V^{(\ell)}_{*}-\m V^{(\ell)}_t \rangle + \frac{ \varrho^\ell}{2 }\| \m V^{(\ell)}_{*}-\m V^{(\ell)}_t\|_F^2, 
	\end{eqnarray}
	Choosing $t = t_s-1$ in the above inequality and letting $s$ goes to $\infty$, we obtain 
	\begin{eqnarray}\label{eq:sup}
	\limsup_{s\rightarrow\infty}  J_{2,\ell}(\m V^{(\ell)}_{t_s}) \leq   J_{2,\ell}(\m V^{(\ell)}_*) .
	\end{eqnarray}
	Here, we have used the fact that $ \nabla_{\m V^{(\ell)}} \bar{\Lc}$ is a gradient Lipchitz continuous function w.r.t. $\m V^{(\ell)}$, the sequence  $\m V^{(\ell)}_{t}$ is bounded and that the distance between two successive iterates tends to zero. Now, we combine \eqref{eq:inf} and \eqref{eq:sup} to obtain
	\begin{eqnarray}\label{eq:limu}
	\lim_{s\rightarrow\infty}  J_{2,\ell}(\m V^{(\ell)}_{t_s}) =  J_{2,\ell}(\m V^{(\ell)}_*)~~ \text{for all}~~ \ell=1, \cdots, K.
	\end{eqnarray}
	Arguing similarly with other variables,  we obtain 
	\begin{subequations}
		\begin{eqnarray}\label{eq:limfull}
		\lim_{s\rightarrow\infty}  J_{1}(\T{G}_{k_s}) &=&  J_{1}(\T{G}_*), \label{eq:limfull1} \\
		\lim_{s\rightarrow\infty}  \bar{\Lc} (\Upsilon_{t_s}) &=&  \bar{\Lc} (\Upsilon^*), \label{eq:limfull4}
		\end{eqnarray} 
	\end{subequations}
	where \eqref{eq:limfull1} follows since $ J_{1}$ is lower semi-continuous; \eqref{eq:limfull4} is obtained from the continuity of function $ \bar{\Lc}$. Thus, $\lim_{s\rightarrow\infty} \Lc ( \Upsilon_{t_s})=\Lc ( \Upsilon_*)$.

	Next, we show that $\Upsilon_*$ is a critical point of $\Lc (.)$. By the first-order optimality condition for the augmented Lagrangian function, we have
	\begin{eqnarray} \label{a55}
	\nonumber
	&& \partial J_{2,\ell}(\m V^{(\ell)}_{t+1}) + \nabla_{\m V^{(\ell)}} \bar{\Lc}(\Upsilon^{\m V^{(\ell)}_{t+1}}_{t+1}) \in \partial_{\m V^{(\ell)}} \Lc(\m V^{(\ell)}_{t+1}),\qquad \ell=1,\dots, K,\\
	\nonumber
	&&\partial J_{1}(\T{G}_{t+1}) + \nabla_{\T{G}} \bar{\Lc}(\Upsilon^{\T{G}_{t+1}}_{t+1}) \in \partial_{\T{G}} \Lc(\T{G}_{t+1}),\\
	&& \gamma  \big({\T{G}}_{t+1} \times_1 {\m V}^{(1)}_{t+1} \cdots \times_K {\m V}^{(K)}_{t+1} -\T{Z}_{t+1} \big) = -  \nabla_{\T{Y}} \Lc(\T{Y}^{t+1}).
	\end{eqnarray}
	Similarly, by the first-order optimality condition for  subproblems \eqref{eq:finalup1}--\eqref{eq:finalup3}, we have
	\begin{eqnarray} \label{a56}
	\nonumber
	&&   \partial J_{2,\ell}(\m V^{(\ell)}_{t+1}) + \nabla_{\m V^{(\ell)}} \bar{\Lc}(\Upsilon^{\m V^{(\ell)}_t}_t)  + \rho^{(\ell)} (\m V_t^{(\ell)}-\m V_{t+1}^{(\ell)}) =0,  \qquad \ell=1, \ldots, K,\\
	&& \partial J_{1}(\T{G}_{t+1}) + \nabla_{\T{G}} \bar{\Lc}(\Upsilon^{\T{G}_{t}}_{t}) + \rho^{\T{G}} (\T{G}_t-\T{G}_{t+1}) =0.
	\end{eqnarray}	
	Combine \eqref{a55} with \eqref{a56} to obtain
	\begin{equation}\label{a59}
	(\eta_{t+1}^1, \dots, \eta_{t+1}^K, \eta_{t+1}^{\T{G}},\eta_{t+1}^{\T{Z}}, \eta_{t+1}^{\T{Y}})\in \partial \Lc (\Upsilon_{t+1}),
	\end{equation}
	where
	\begin{eqnarray} \label{a57}
	\nonumber
	&& \eta^{t+1}_{\m V^{(\ell)}} :=  \nabla_{\m V^{(\ell)}} \bar{\Lc}(\Upsilon^{\m V^{(\ell)}_{t+1}}_{t+1}) - \nabla_{\m V^{(\ell)}} \bar{\Lc}(\Upsilon^{\m V^{(\ell)}_t}_t)  - \rho^{(\ell)} (\m V_t^{(\ell)}-\m V_{t+1}^{(\ell)}),  \qquad \ell=1, \ldots, K, \\
	\nonumber
	&& \eta^{t+1}_{\T{G}} := \nabla_{\T{G}} \bar{\Lc}(\Upsilon^{\T{G}_{t+1}}_{t+1}) - \nabla_{\T{G}} \bar{\Lc}(\Upsilon^{\T{G}_{t}}_{t}) - \rho^{\T{G}} (\T{G}_t-\T{G}_{t+1}),\\
	\nonumber
	&& \eta^{t+1}_{\T{Z}} := \nabla_{\T{Z}} \bar{\Lc}(\Upsilon^{\T{Z}_{t+1}}_{t+1}) = \T{Y}_t-\T{Y}_{t+1},\\
	&&  \eta^{t+1}_{\T{Y}} := \frac{1}{\gamma}  (\T{Y}_t-\T{Y}_{t+1}).
	\end{eqnarray}
	Note that the function $\bar{\Lc} (\Upsilon)$  is a gradient Lipchitz continuous function w.r.t. $\m V^{(1)} , \dots, \m V^{(K)},\T{G}$. Thus, 
	\begin{eqnarray} \label{aa57}
	\nonumber
	\|\nabla_{\m V^{(\ell)}} \bar{\Lc}(\Upsilon^{\m V^{(\ell)}_{t+1}}_{t+1}) - \nabla_{\m V^{(\ell)}} \bar{\Lc}(\Upsilon^{\m V^{(\ell)}_t}_t) \| & \leq &  \rho^{(\ell)} \|\m V_t^{(\ell)}-\m V_{t+1}^{(\ell)}\|,  \qquad \ell=1, \ldots, K, \\
	\nonumber\\
	\| \nabla_{\T{G}} \bar{\Lc}(\Upsilon^{\T{G}_{t+1}}_{t+1}) - \nabla_{\T{G}} \bar{\Lc}(\Upsilon^{\T{G}_{t}}_{t}) \| & \leq &  \rho^{\T{G}} \| \T{G}_t-\T{G}_{t+1}\|,
	\end{eqnarray}  
	Using \eqref{aa57} and \eqref{a57}, we obtain
	\begin{eqnarray} \label{a66}
	\lim_{t\rightarrow\infty} \big(\|\eta_{t+1}^1\|, \dots, \|\eta_{t+1}^K\|, \|\eta_{t+1}^{\T{G}}\|, \|\eta_{t+1}^{\T{Z}}\|, \|\eta_{t+1}^{\T{Y}}\|\big) =\big(0,\dots, 0\big).
	\end{eqnarray}
	Now, from \eqref{a59} and \eqref{a66}, we conclude that $(0,\dots,0)\in \partial \Lc(\Upsilon_*)$ due to the closure property of $\partial \Lc$. Therefore,  $\Upsilon_*$ is a critical point of $\Lc (.)$. This completes the proof of (i).
	\\
	(ii). The proof follows from \citep[Lemma 5 and Remark 5]{bolte2014proximal}.
	\\
	(iii). Choose $\Upsilon_*\in \T C(\Upsilon_0)$. There exists a subsequence $\Upsilon_{t_s}$ converging to $\Upsilon_{*}$ as $s$ goes to infinity. Since we have proven that $
	\lim\limits_{s\rightarrow\infty}\Lc(\Upsilon_{t_s}) = \Lc(\Upsilon_{*})$, and $\Lc(\Upsilon_t)$ is a non-increasing sequence, we conclude that $\Lc(\Upsilon_*) = \lim\limits_{t\rightarrow\infty} \Lc(\Upsilon_t)=:\T L_*$, hence the restriction of $\Lc(\Upsilon)$ to $\T C(\Upsilon_{0})$ equals to $\T L_*$.
\end{proof}

\subsection{On the Simulated Coupled Datasets}\label{app:c}

	The script used to create coupled tensor matrix out of initial factors is called \texttt{create\_{coupled}} \citep{acar2014structure}. In addition to choosing the distribution of the random simulation of data, one may indicate the size of the tensor as well as the mode to which the matrix is attached.
\begin{figure}[h]
	\centering
 \includegraphics[width=5in]{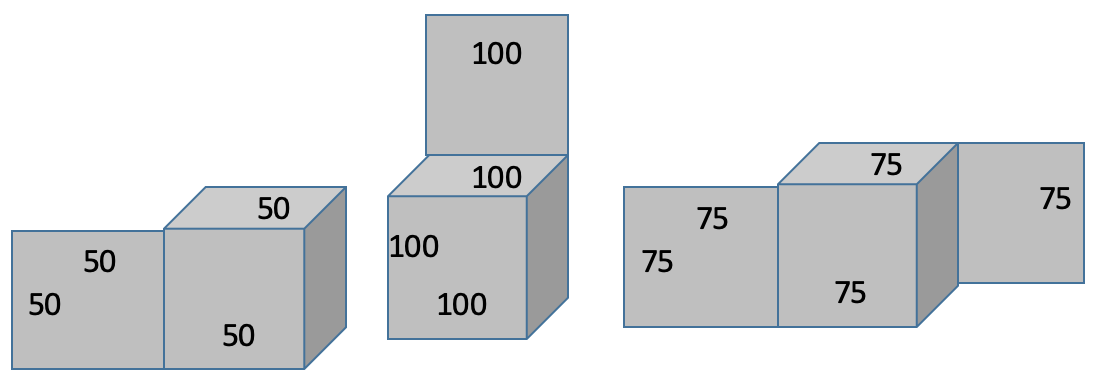}
	\caption{ Using the script  \texttt{create\_{coupled}}  to create simulated coupled tensor matrix datastes from initial factors. Left, the size is $[50~50~50~50]$ and the mode reads $\{[ 1~2~3], [1,4]\}$. Similarly, for the middle figure, the size is $[100~100~100~100]$ while the mode input is $\{[ 1~2~3], [2,4]\}$. Lastly, the figure in right is created by using $[75~75~75~75~75]$ and  $\{[ 1~2~3], [1,4], [3,5]\}$ for the size and mode inputs, respectively. For more details, we refer the reader to \citep{acar2014structure}.  }
	\label{fig:coup}
\end{figure} 
\subsection{On Bilevel Optimization}\label{app:bl}
Bilevel optimization has increased in popularity in recent years. Bilevel problems are often used to model hierarchical processes in which
an upper level problem makes a decision anticipating the rational inputs of the lower level problem. Consider a linear bilevel problem of the following form:
\begin{equation}\label{prob:1}
\begin{aligned}
\min\limits_{x,y} \quad & c^\top x+d^\top y \\
\textrm{s.t.} \quad & \mathbf Ax+\mathbf By\ge a\\
\quad & y\in \argmax\limits_{\bar y}\{e^\top \bar y: \mathbf Cx+\mathbf D\bar y\le b\}, 
\end{aligned}
\end{equation}
We refer to $x$ and $y$ as the upper and lower-level variables, respectively. Here, 
\begin{align*}
&a\in \mathbb R^k, \quad b\in \mathbb R^\ell \quad c\in \mathbb R^n \quad d,e\in \mathbb R^m\\
& \mathbf A\in \mathbb R^{k\times n}\quad 
\mathbf B\in \mathbb R^{k\times m}\quad
\mathbf C\in \mathbb R^{\ell\times n}\quad
\mathbf D\in \mathbb R^{\ell\times m}
\end{align*}
Considered as a parametric linear optimization problem (LP), the lower-level
problem is
\begin{equation}
\begin{aligned}
\max\limits_{y} \quad & e^\top y\\
\textrm{s.t.} \quad & \mathbf D y\le b-\mathbf C\bar x, 
\end{aligned}
\end{equation}
and the dual problem is given by:
\begin{equation}
\begin{aligned}
\min\limits_{\lambda} \quad & (b-\mathbf C\bar x)^\top \lambda\\
\textrm{s.t.} \quad & \mathbf D^\top \lambda =e
\qquad \&\qquad  \lambda \ge 0.
\end{aligned}
\end{equation}
Note that the dual polyhedron $D := \{ \lambda \in \mathbb R^\ell,\mathbf D^\top \lambda =e, \lambda \ge 0\}$  does not depend
on primal upper-level variables.
The Strong Dualtity Theorem~\ref{thm:sdt} states a key result about the relationship between bilevel and single-level linear optimization problems:
	\begin{theorem}[Strong Duality Theorem ]\label{thm:sdt}
	\citep{dantzig2006linear}
		Considering the linear optimization (LP) problem given in Eq.~\eqref{prob:1}, there are four possibilities:
		\begin{enumerate}
			\item Both primal and dual have no feasible solutions (are infeasible).
			\item The primal is infeasible and the dual unbounded.
			\item The dual is infeasible and the primal unbounded.
			\item Both primal and dual have feasible solutions and their values are equal.
		\end{enumerate}
	\end{theorem} Applying this theorem to the lower-level problem, one can \emph{equivalently} reformulate the bilevel problem \ref{prob:1} to the single-level problem
\begin{equation}\label{prob:2}
\begin{aligned}
\min\limits_{x,y,\lambda} \quad & c^\top x+d^\top y \\
\textrm{s.t.} \quad & \mathbf Ax+\mathbf By\ge a\quad \& \quad  \mathbf Cx+\mathbf Dy\le b \quad \& \quad \mathbf D^\top \lambda=e\\
\quad &  \lambda\ge 0\quad \& \quad e^\top y -b^\top \lambda +x^\top \mathbf C^\top \lambda\ge 0.
\end{aligned}
\end{equation}
Problems \ref{prob:1} and \ref{prob:2} are equivalent. 
For a general ADMM problem,  we consider an optimization
problem in the specific form
\begin{equation}\label{prob:admm}
\begin{aligned}
\min\limits_{\mathbf x,\mathbf y} \quad & f(\mathbf x,\mathbf y)\\
\textrm{s.t.} \quad & g(\mathbf x,\mathbf y)=0, h(x,y)\ge 0 \\
\quad & \mathbf x\in X, \quad \mathbf y\in Y
\end{aligned}
\end{equation}
where $\mathbf x\in\mathbb R^n$, $\mathbf y\in\mathbb R^m$ are variable vectors. The feasible set of this problem is:
$$ \Omega:= \{(\mathbf x,\mathbf y)\in X\times Y: g(\mathbf x, \mathbf y)=0, h(\mathbf x, \mathbf y)\ge 0\},$$
For discussing the theoretical properties of ADMMs, we need the following assumption:
\begin{remark}[assumptions]
	The objective function $f:\mathbb R^n\times \mathbb R^m\to \mathbb R$ and the constraint $g:\mathbb R^n\times \mathbb R^m\to \mathbb R^k$ and $h:\mathbb R^n\times \mathbb R^m\to \mathbb R^\ell$ are continuous and the sets $X$
	and $Y$ are non-empty and compact.
\end{remark}
The following general convergence result holds:
\begin{theorem}
	Let $\{(\mathbf x^i, \mathbf y^i)\}_{i=0}^\infty$ be a sequence with $(\mathbf x^i, \mathbf y^i)\in \Sigma(\mathbf x^i, \mathbf y^i)$, where \begin{multline}
	\Sigma(\mathbf x^i, \mathbf y^i):=
	\{ (\mathbf x^*, \mathbf y^*): f(\mathbf x^*, \mathbf y^i)\le f(\mathbf x, \mathbf y^i), \forall \mathbf x\in X; f(\mathbf x^*, \mathbf y^*)\le f(\mathbf x^*, \mathbf y), \forall \mathbf y\in Y\}
	\end{multline}
	Suppose that Assumption 1 holds and that the solution of the first optimization
	problem is always unique. Then, every convergent subsequence of $\{(\mathbf x^i, \mathbf y^i)\}_{i=0}^\infty$ converges to a partial minimum. For two limit points $z, z'$ of such subsequences it
	holds that $f(z) = f(z')$.
\end{theorem}
Stronger convergence results can be obtained
if stronger assumptions on $f$ and $\Sigma$ are made. 
\begin{remark}[Corollary]
	Suppose that the assumptions of Theorem 1 are satisfied. Then, the
	following holds:
	\begin{enumerate}
		\item If $f$ is continuously differentiable, then every convergent subsequence of $\{(\mathbf x^i, \mathbf y^i)\}_{i=0}^\infty$ converges to a stationary point of Problem \ref{prob:admm}. 
		\item If $f$ is continuously differentiable and if $f$ and $\Sigma$ are convex, then every
		convergent subsequence of $\{(\mathbf x^i, \mathbf y^i)\}_{i=0}^\infty$ converges to a global minimum of
		Problem \ref{prob:admm}. 
	\end{enumerate}
\end{remark}
Problem \ref{prob:admm} can be seen as a quasi block-separable
problem, where the blocks are given by the variables $\mathbf x$ and $\mathbf y$ as well as their
respective feasible sets $X$ and $Y$. Here, quasi means that there still
are the constraints $g$ and $h$ that couple the feasible sets of the two blocks. The main
idea of ADMM is to alternately solve in the directions of the blocks separately. Now, further relaxing the coupling constraints, $g$ and $h$, with the following penalty function
\begin{equation}
\varphi(\mathbf x, \mathbf y; \mu, \nu):= f(\mathbf x, \mathbf y)+\sum\limits_{t=1}^k \mu_k|g_t(\mathbf x, \m y)|+ \sum\limits_{t'=1}^\ell \nu_{t'}[h_{t'}(\mathbf x, \m y)]^-, 
\end{equation}
where $\mu$ and $\nu$ are vector penalty parameters of sizes $k$ and $\ell$, respectively. Moreover, $[\cdot]^-:=\max\{0,\cdot\}$. 
The penalty ADMM consists of an inner and an outer
loop. Applying the ADMM algorithm to the inner problem, if the inner loop iteration terminates with a partial minimum, if the coupling constraints are satisfied, this terminates the process.
\begin{theorem}
	Suppose that Assumption 1 holds, and that $\mu^i_t$, $t=1:k$ and ${\nu^i}_{t'}$, $t'=1:\ell$ are monotonically increasing sequences. Let $\{(\mathbf x^i, \mathbf y^i)\}_{i=0}^\infty$ be a sequence partial minima of the following problem
	\begin{equation}
	\min\limits_{\mathbf x, \mathbf y}\varphi(\mathbf x, \mathbf y; \mu, \nu), \quad \text{s.t.}\quad \mathbf x\in X, \quad \mathbf y\in Y. 
	\end{equation}
	Then, there exist weights $\bar \mu, \bar \nu\ge 0$, such that $(\mathbf x^*, \mathbf y^*)$ is a partial minimizer of the weighted feasibility measure, if is feasible to Problem \ref{prob:admm}
	\begin{itemize}
		\item if $f$ is continuous, then $(\mathbf x^*, \mathbf y^*)$ is a partial minimum of Problem \ref{prob:admm}, 
		\item if $f$ is continuously differentiable, then $(\mathbf x^*, \mathbf y^*)$ is a stationary point of Problem \ref{prob:admm}, 
		\item if $f$ is continuously differentiable, and $f$ and $\Sigma$ are convex, then $(\mathbf x^*, \mathbf y^*)$ is a global optimum of Problem \ref{prob:admm}. 
	\end{itemize}
\end{theorem}
\end{document}